\documentclass[draftcls, onecolumn]{IEEEtran}
\ifCLASSINFOpdf
\else
\fi

\usepackage{amsmath,amsthm}
\usepackage{amsfonts}

\usepackage{graphicx}
\usepackage{subfigure}
\usepackage{float}

\graphicspath{{Codes/}}

\newcommand{\Cbb} {\mathbb{C}}
\newcommand{\R} {\mathbb{R}}

\newcommand{\Abs}[1] {\left| #1 \right|}

\newtheorem{thm}{Theorem}[section]
\newtheorem{cor}[thm]{Corollary}
\newtheorem{lem}[thm]{Lemma}
\newtheorem{prop}[thm]{Proposition}
\theoremstyle{definition}
\newtheorem{defn}[thm]{Definition}
\theoremstyle{remark}
\newtheorem{rem}[thm]{Remark}
\numberwithin{equation}{section}

\begin{document}
%
\title{Asymptotics of empirical eigenvalues for large separable covariance matrices}
%
%
%

\author{Tiebin~Mi,~\IEEEmembership{Member,~IEEE,}
        Robert~Caiming~Qiu,~\IEEEmembership{Fellow,~IEEE,}
\thanks{T. Mi is with the School of Electronics, Information and Electrical Engineering (SEIEE), Shanghai Jiaotong University, 200240, China (e-mail: mitiebin@sjtu.edu.cn).}%
\thanks{R. Qiu is with the School of Electronics, Information and Electrical Engineering (SEIEE), Shanghai Jiaotong University, 200240, China (e-mail: rcqiu@sjtu.edu.cn).}%
}

\maketitle

\begin{abstract}
We investigate the asymptotics of eigenvalues of sample covariance matrices associated with a class of non-independent Gaussian processes (separable and temporally stationary) under the Kolmogorov asymptotic regime. The limiting spectral distribution (LSD) is shown to depend explicitly on the Kolmogorov constant (a fixed limiting ratio of the sample size to the dimensionality) and parameters representing the spatio- and temporal- correlations. The Cauchy, M- and N-transforms from free harmonic analysis play key roles to this LSD calculation problem. The free multiplication law of free random variables is employed to give a semi-closed-form expression (only the final step is numerical based) of the LSD for the spatio-covariance matrix being a diagonally dominant Wigner matrix and temporal-covariance matrix an exponential off-diagonal decay (Toeplitz) matrix. Furthermore, we also derive a nonlinear shrinkage estimator for the top eigenvalues associated with a low rank matrix (Hermitian) from its noisy measurements. Numerical studies about the effectiveness of the estimator are also presented.
\end{abstract}


%
\IEEEpeerreviewmaketitle

\section{Introduction}\label{S:Introduction}
Many applications in signal processing and other sciences involve modeling a dependency across a wide range of spatial and temporal scales, e.g., \cite{vuran2004spatio, bolcskei2006space, xu2002spatial, gneiting2006geostatistical, barndorff2004levy}. Traditionally, statistical theory has been concerned with studying the properties of a spatio-temporal process as the sample size tends to \emph{infinity}, while the number of variables is fixed. Nevertheless, there are always a few snapshots available to us about the underlying process in practice. A key challenge in spatio-temporal data analysis is therefore to summarize or infer the process’ behavior with these \emph{finite} samples. In this article, we will discuss the limiting behavior of the empirical spectral distribution associated with a class of Gaussian processes (large, non-independent, separable and temporally stationary).

For convenience, we denote by $U(s, t)$ a spatio-temporal process. Particularly, the covariance function $\text{cov}(U(s, t), U(s', t'))$ provides us a succinct but informative summary on the process. If no further structure exists, $\text{cov}(U(s, t), U(s', t'))$ depends entirely on the space-time indexes $(s, t)$ and $(s', t')$. However, stationary and separability are often imposed to ensure that the covariance function could be manipulated easily \cite{cressie2015statistics}.

Specifically, a spatio-temporal random function possesses temporally stationary covariance if $\text{cov}(U(s, t), U(s', t'))$ depends on the time index only through the temporal lag $t-t'$ \cite{gneiting2006geostatistical}. Similarly, it has spatially stationary covariance if $\text{cov}(U(s, t), U(s', t'))$ depends on $s$ and $s'$ only through the spatial separation $s-s'$. A spatio-temporal process has stationary covariance if it has both spatially and temporally stationary covariance. Under this assumption, there exists a function $c(\cdot, \cdot)$ such that $\text{cov}(U(s, t), U(s', t')) = c(s-s', t-t')$.

The separability, which offers significantly computational benefits, is a desirable property for the spatio-temporal models \cite{cressie2015statistics, fuentes2006testing, sherman2011spatial}. A process $U(s, t)$ has a separable covariance function if $\text{cov}(U(s, t), U(s', z')) = \sigma_{s} (s, s') \sigma_{t} (t, t')$ for some spatio- and temporal- covariances $\sigma_{s}(\cdot, \cdot)$ and $\sigma_{t}(\cdot, \cdot)$. Roughly speaking, with the separability, the spatio- and temporal-correlations are ``decoupled'' completely. Moreover, with both the separability and temporal stationary, the associated covariance function could be written as $\text{cov}(U(s, t), U(s', z')) = \sigma_{s} (s,s') \sigma_{t} (t-t')$.

In this article we will employ a framework inspired by random matrix theory to characterize the asymptotic behaviour of \emph{large sample covariance matrices} associated with a class of separable and temporally stationary Gaussian processes (noises). The setting in random matrix theory differs substantially from traditional situation where the number of variables is fixed, but the sample size for examination tends to infinity. The proposed asymptotic approach allows us to obtain the limiting of the eigenvalues (under the condition that the number of variables and samples size both tend to infinity) and predict the behavior of the empirical eigenvalues associated with large sample covariance matrices with reasonable accuracy. This is the primary motivation for this article.

We begin with the definition of separable Gaussian process. Specifically, for the samples of a separable Gaussian process, we consider a matrix $U_{N,T}$ with the form $U_{N,T} = \left( \Sigma^{s}_{N} \right)^{\frac{1}{2}} W_{N,T} \left( \Sigma^{t}_{T} \right)^{\frac{1}{2}}$, where $W_{N,T}$ is a random matrix consisting of $N \times T$ independent identically distributed (i.i.d.) Gaussian random variables, $\Sigma^{s}_{N}$ and $\Sigma^{t}_{T}$ are $N \times N$ and $T \times T$ non-negative definite matrices, representing spatio- and temporal- covariances, respectively. Clearly, with the impact of $\left( \Sigma^{s}_{N} \right)^{\frac{1}{2}}$ and $\left( \Sigma^{t}_{T} \right)^{\frac{1}{2}}$, $U_{N,T}$ consists of non-independent Gaussian random variables. For the associated sample covariance matrix $C_N$, we mean $C_N = \frac{1}{T} U_{N,T} U_{N,T}^T = \frac{1}{T} \left( \Sigma^{s}_{N} \right)^{\frac{1}{2}} W_{N,T} \Sigma^{t}_{T} W_{N,T}^{\text{T}} \left( \Sigma^{s}_{N} \right)^{\frac{1}{2}}$.

To examine the asymptotic statistics of large sample covariance matrices, a new approach based on free probability theory (see for instance \cite{burda2010random, burda2011applying}) is developed under the Kolmogorov condition, i.e., $N, T \to \infty$ and $\frac{N}{T} \to c$ (Kolmogorov constant). We emphasize that this approach is different greatly from the methods in \cite{paul2009no, couillet2014analysis}. The Cauchy, M- and N-transforms from free harmonic analysis play key roles in this spectral distribution estimation problem. We will investigate how to calculate with these transforms a semi-closed-form expression of the LSD of the sample covariance matrices associated with a class of separable Gaussian processes.

Before we proceed, let us first discuss briefly the computational and statistical challenges we are faced under the Kolmogorov condition.

\subsection{Challenges and opportunities from massive data}
In recent years, modern signal processing techniques often amplify system's pressure from dealing with massive data in pursuit of high performance. Massive data are always featured in high-dimensionality and large sample size, i.e., the dimensionality of variables is comparable to the size of available samples (the Kolmogorov condition falls into this regime also). That will introduce computational and statistical challenges, including spurious correlation, noise accumulation and scalability.

A good illustration for massive data processing is the development in wireless communication \cite{bi2015wireless, rusek2013scaling}. Basically, the more antennas the transmitter/receiver are equipped with, the more degrees of freedom could be achieved. The LTE standard allows for up to eight antennas at the base station. The massive MIMO goes even further, which consists of a very large number of antennas (hundreds or thousands), simultaneously serves a great number of user terminals (tens or hundreds). Therefore, it's necessary to know what will happen for a large system limit, where the numbers of antennas and user terminals grow infinitely large with the same speed \cite{hoydis2011massive, rusek2013scaling}.

Unfortunately, the Kolmogorov condition limits the applicability of many well-known approaches \cite{fan2014challenges}. We consider the covariance matrix estimation as an example. Given an $N \times T$ data matrix $X$, consisting of $T$ samples drawn from multivariate Gaussian distribution $\mathcal{N}(0, I_{N})$, it is well-known that the sample covariance matrix is a good approximate of the population covariance matrix when $N$ remains fixed and $T \to \infty$, i.e., $\lim_{T \to \infty} \frac{1}{T} X X^* = I_{N}$.

By contrast, for the situation where the dimension $N$ is comparable to the sample size $T$, many new phenomena arise. Perhaps the first high-dimensionality phenomenon is that, if $N \simeq T$, the empirical eigenvalues do \emph{not} converge to their population counterparts even as the sample size grows. The famous Mar\u{c}enko-Pastur law (illustrated in Theorem~\ref{T:MP_Law}) gives a celebrated function for the asymptotic behavior of eigenvalues of a large sample covariance matrix associated with multivariate Gaussian distributions. We can see from Fig.~\ref{F:MP_Law} that the eigenvalues of the sample covariance matrices are more spread out, or dispersed, than the population eigenvalues.

Such a paradigm shift to big/massive data in signal processing will lead to significant progresses on the development of new theory and method. Although there exist various distinct statistical problems, we primarily focus on the limiting spectral statistics associated with large sample covariance matrices. To this end, a framework inspired by random matrix theory and several transforms from free harmonic analysis will be adapted.

The organization of this article is as follows. A brief description of random matrix theory is given in Section~\ref{S:RMT}. The roles of Cauchy transform in the calculation of Mar\v{c}enko-Pastur distribution and top eigenvalues of a low-rank spiked covariance model are also discussed. In Section~\ref{S:FP}, we present some preliminaries for free probability theory. For those who are familiar with this may skip most of the section. Section~\ref{S:Separability} focuses on the models of spatio- and temporal- covariance matrices. We restrict the spatio- and temporal-covariance matrices to a diagonally dominant Wigner matrix and a bi-infinite Toeplitz matrix, respectively. Section~\ref{S:Cauchy_M_N} is dedicated to the calculation of Cauchy, M- and N-transforms associated with the Mar\v{c}enko-Pastur distribution, diagonally dominant Wigner matrix and bi-infinite Toeplitz matrix. All of the associated Cauchy, M- and N-transforms are analytical. Finally, the article turns to the calculation of LSD of large sample covariance matrices associated with autoregressive Gaussian processes. We also have a very short discussion on the nonlinear shrinkage estimator for the top eigenvalues of a low rank matrix (Hermitian) from its noisy measurements in Section~\ref{S:LSD_AR}.

\section{Random matrix theory and Cauchy transform}\label{S:RMT}
Random matrix theory, which studies how the eigenvalues behave asymptotically when the order of an underlying random matrix grows, is gaining increasing attention for analyzing complex high-dimensional data. Early interest in random matrices arose in the context of multivariate statistics, it was Wigner who introduced random matrix ensembles and derived the first asymptotic result \cite{wigner1958distribution, akemann2011oxford}. Since then, there have been a large number of literatures studying the asymptotic behavior of various random matrices. In this section, we will introduce some necessary concepts in random matrix theory. For a thorough treatment on random matrix theory, we refer the reader to \cite{el2008spectrum, bai2010spectral, tao2012topics, yao2015sample, couillet2011random, anderson2010introduction, tulino2004random, mingo2017free, speicher2014free}.

\begin{defn}\label{D:ESD}
Given an $N \times N$ Hermitian matrix $M_N$, the empirical spectral distribution (or ESD for short) is defined as, for $x \in \R$,
\[
\mu_{M_N}(x) \equiv \frac{1}{N} \sum_{i=1}^{N} 1_{ \lambda_i(M_N) \le x }
\]
where $\lambda_1(M_N) \le \cdots \le \lambda_N(M_N)$ are the (necessarily real) eigenvalues of $M_N$, including multiplicity.
\end{defn}

Note that $\mu_{M_N}(x)$ is actually a probability measure on $\R$. It is non-decreasing with $\lim_{x \to -\infty} \mu_{M_N}(x) = 0$ and $\lim_{x \to \infty} \mu_{M_N}(x) = 1$.

Random matrix theory asserts that there exist some classes of Hermitian random matrices $M_N$ whose (random) ESD $\mu_{M_N}(x)$ converges, as $N \to \infty$, towards a non-random function $\mu(x)$. This function $\mu$, if it exists, is called the limiting spectral distribution (LSD for short). Moreover, for many Hermitian random matrices, the corresponding LSD is supported on a bounded interval (bulk) and illustrating sharp edges on the boundary. The eigenvalue significantly outside of the interval (noise) will be a good indicator of the non-random behavior (signal).

We now introduce Cauchy transform, which uses complex-analytic methods and has turned out to be a powerful tool in the calculation of LSDs.

\begin{defn}
Let $\mathbb{C}^+ = \{ z \in \mathbb{C} | \text{Im} (z) > 0 \}$ be the complex upper half-plane and $\mu$ a probability measure on $\R$. The Cauchy transform $\mathcal{G}_{\mu}(z)$ is an analytic function defined as, for $z \in \mathbb{C}^+ \cup \mathbb{R}\backslash \text{supp} (\mu)$,
\[
\mathcal{G}_{\mu}(z) \equiv \int_{\R} \frac{1}{z-t} d \mu(t).
\]
\end{defn}

A probability measure could be retrieved from its Cauchy transform via the so-called Stieltjes inversion formula. Particularly, if $x$ is a continuous point of $\mu$,
\begin{equation}\label{E:StieltjesInversion}
\mu(x) = - \frac{1}{\pi} \lim_{y \to 0^+} \int_{-\infty}^{x} \text{Im} \left( \mathcal{G}_{\mu} (x + iy) \right) dx.
\end{equation}
With this Stieltjes inversion formula, the statistical information on a probability measure is completely encoded in the corresponding Cauchy transform function. However, it is worth noting that, for an arbitrary probability distribution, it might not be possible to obtain an analytical or even an easy numerical solution via the Stieltjes inversion formula.

Next we will consider Cauchy transform associated with a Hermitian random matrix. By Definition~\ref{D:ESD}, we know that $d \mu_{M_N}(t)$ is the summation of a series of Dirac delta functions that concentrate around the eigenvalues of $M_N$. The associated Cauchy transform is then defined as, for $z \in \mathbb{C}^+ \cup \mathbb{R}\backslash \text{supp} (\mu_{M_N})$,
\begin{equation*}
\begin{aligned}
\mathcal{G}_{M_N} (z) & \equiv \int_{\R} \frac{1}{z-t} d \mu_{M_N}(t) \\
& = \text{tr} \left( (z I_N - M_N)^{-1} \right) = \sum_{k \ge 0} \frac{\text{tr} \left( M_N^k \right) }{z^{k+1}} ,
\end{aligned}
\end{equation*}
where $\text{tr}(M_N) \equiv \frac{1}{N} \text{Tr}(M_N) $ is the normalized trace.

\subsection{Mar\v{c}enko-Pastur law}
The Cauchy transform was used by Mar\v{c}enko and Pastur to derive the well-known Mar\v{c}enko-Pastur law of large dimensional Gram matrices of random matrices with i.i.d. Gaussian entries \cite{marvcenko1967distribution}.

\begin{thm}[\cite{marvcenko1967distribution}]\label{T:MP_Law}
Consider an $N \times T$ Gaussian random matrix $W_{N \times T}$ (consists of i.i.d. Gaussian entries of mean 0 and variance 1). Under the Kolmogorov condition, i.e., $N, T \to \infty$ and $\frac{N}{T} \to c \in (0, \infty)$, the ESD of $C_{N} = \frac{1}{T} W_{N \times T} W_{N \times T}^{\text{T}}$ converges weakly and almost surely to a non-random distribution function (Mar\v{c}enko-Pastur distribution) with density function
\[
d \mu (x) = (1 - \frac{1}{c})^{+} \delta(x) + \frac{1}{2 \pi c x} \sqrt{(x-a)^{+}(b-x)^{+}}
\]
where $a = (1 - \sqrt{c})^2$, $b = (1 + \sqrt{c})^2$, $a^{+} = \max (0, a)$ and $\delta(x) = 1_{\{ 0 \}} (x)$.
\end{thm}

\begin{figure}
\centering
\subfigure[Eigenvalue spreading phenomenon ($N = 4000, T = 8000$).]{
\label{F:MP_Law:1}
\includegraphics[width=0.47\columnwidth]{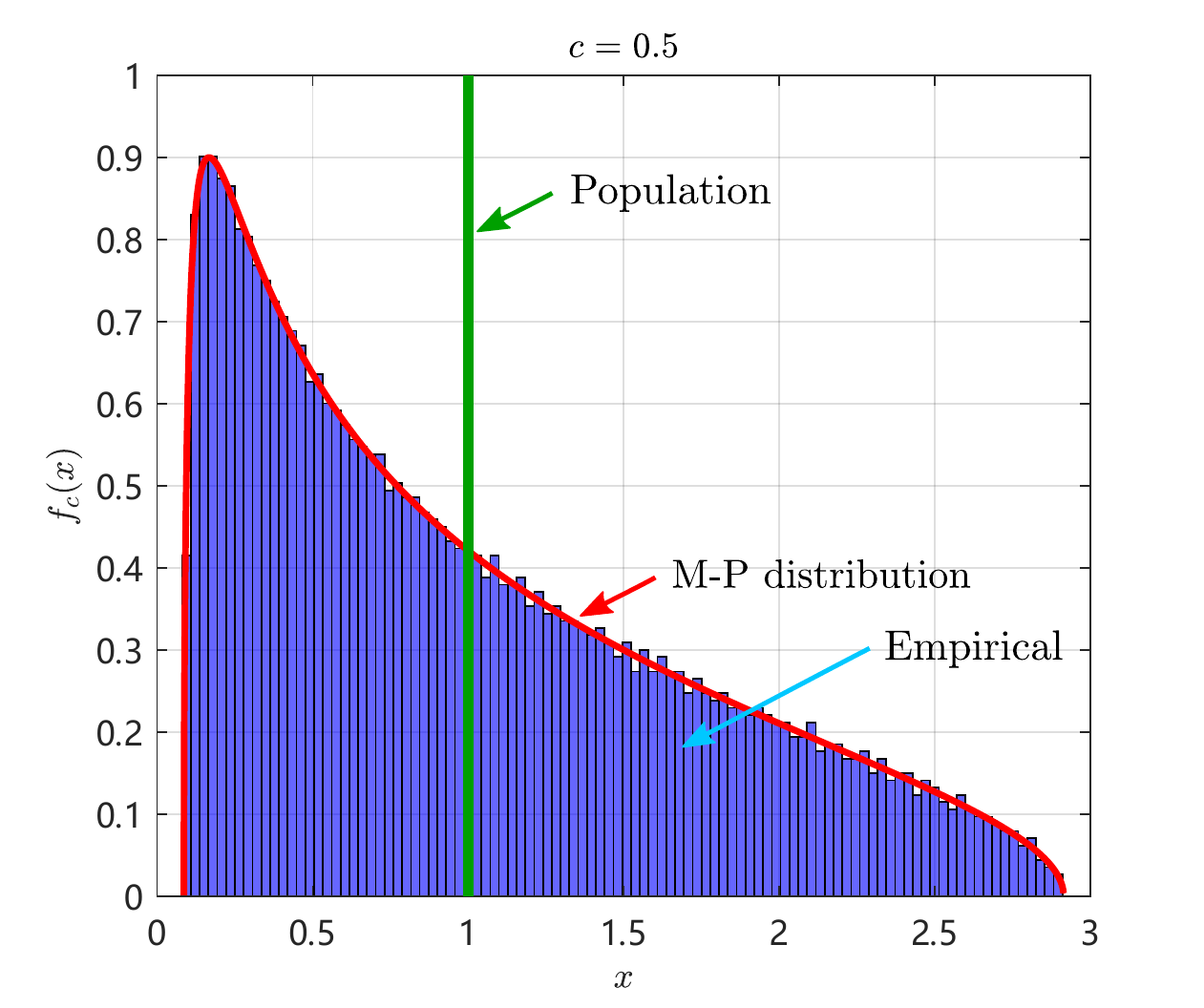}}
\subfigure[Three instances of the Mar\v{c}enko-Pastur distribution for different values of the limiting ratio $c$.]{
\label{F:MP_Law:2}
\includegraphics[width=0.47\columnwidth]{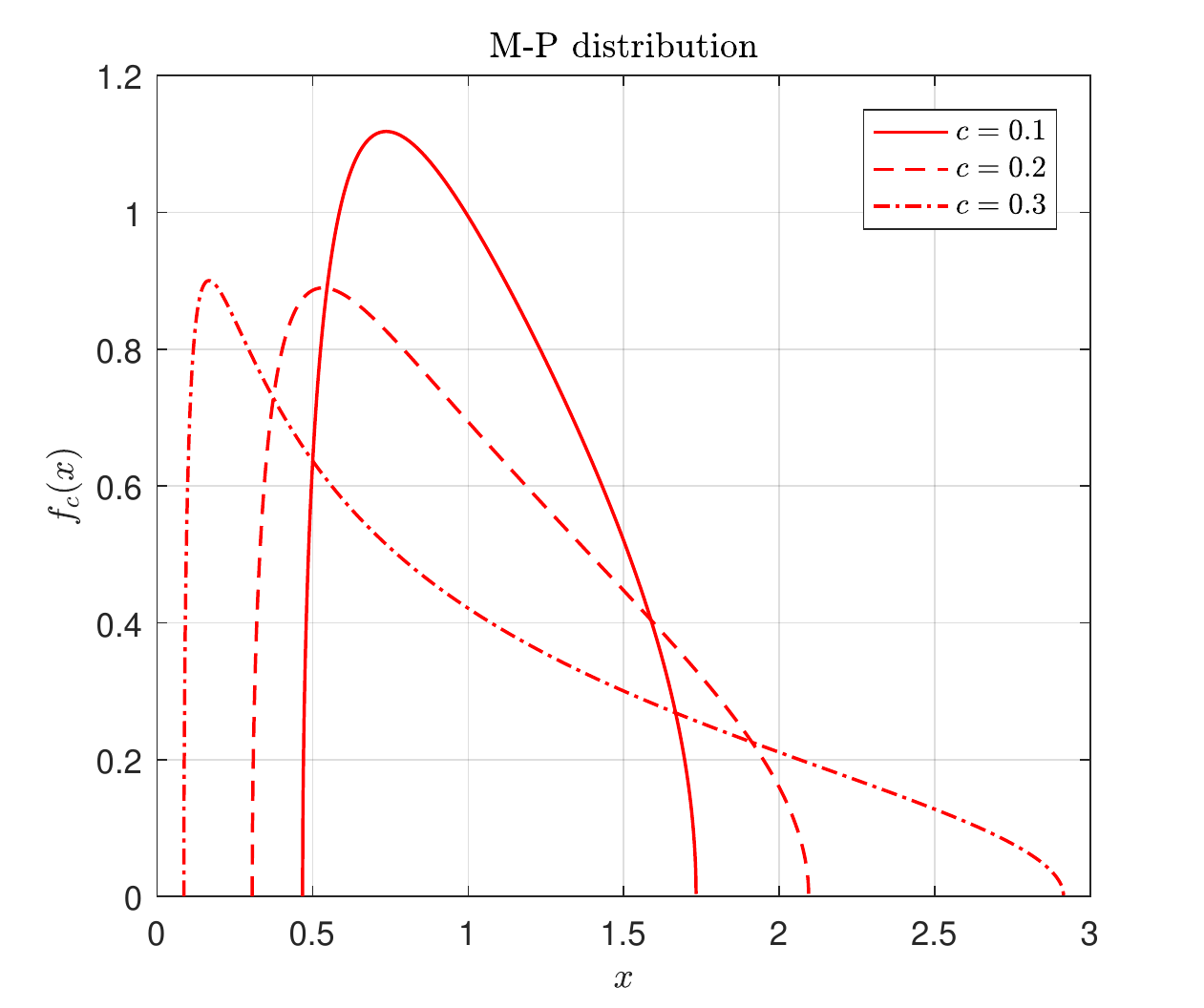}}
\caption{High-dimensional phenomenon: spreading of eigenvalues. Population eigenvalues, LSDs and ESDs are shown in green, red and blue, respectively.}
\label{F:MP_Law}
\end{figure}

From Theorem~\ref{T:MP_Law}, we see that the Mar\v{c}enko-Pastur distribution is supported on $\left[(1 - \sqrt{c})^2, (1 + \sqrt{c})^2 \right]$. Thus one of the high-dimensional phenomena relates to Kolmogorov condition is that the sample eigenvalues are more spread out than the population eigenvalues for large $N$ and $T$. Fig.~\ref{F:MP_Law:1} illustrates this phenomenon for $N = 4000$ and $T = 8000$. All population eigenvalues are equal to 1 (in green), but the 4000 sample eigenvalues have a histogram (in blue) that spreads over a range. Fig.~\ref{F:MP_Law:2} illustrates the Mar\v{c}enko-Pastur law for different values of the Kolmogorov constant $c$.

\subsection{A low-rank spiked sample covariance model}\label{SS:LowRankWigner}
In recent years, some authors investigated in detail the low rank ``spiked'' covariance matrix model, e.g., \cite{johnstone2018pca, johnstone2001distribution, donoho2018optimal, wang2017asymptotics, loubaton2011almost, baik2006eigenvalues, shabalin2013reconstruction, hachem2013subspace, vallet2012improved, bun2017cleaning, biroli2019large, bun2016rotational}. A spiked sample covariance matrix $Y_N$ consists of a known ``base'' covariance matrix $C_N$ (bulk, noise) and a low-rank perturbation $X_N$ (spikes, signals)
\[
Y_N = X_N + C_N = \sum_{i = 1}^{r} \theta_i \phi_{i} \phi_{i}' + C_N,
\]
where the eigenvalues of $X_N$, or the signal strengths $\theta_1, \ldots, \theta_i$ are often taken as unknown. We assume that the rank $r$ is relatively small and stays fixed in the asymptotic regime. In particular, there exist a notable high-dimensional phenomenon - the eigenvalue phase transition.

In \cite{benaych2011eigenvalues}, Benaych-Georges, F. and Nadakuditi, R.R. focused on the role of Cauchy transform in calculating the top eigenvalues of the low-rank spiked covariance model. They showed that the limiting non-random eigenvalues associated with the perturbed matrix depend explicitly on the levels of the spikes and the LSD of the bulk via the Cauchy transform.

\begin{thm}[Eigenvalue phase transition for observation \cite{benaych2011eigenvalues}]\label{T:EigenvalueNonlinearMappingForward}
Consider a sample covariance matrix $Y_N = X_N + C_N = \sum_{i = 1}^{r} \theta_i \phi_{i} \phi_{i}' + C_N$. Without loss of generality, we assume that $r \ll N$, $\theta_1 \ge \ldots \ge \theta_r \ge 0$ and, as $N \to \infty$, $\lambda_{i} (Y_N) \xrightarrow{\text{a.s.}} \eta_i $, $\mu_{C_N} \xrightarrow{\text{a.s.}} \mu_{C}$ (compactly supported). Let $a$ be the supremum of the support of $\mu_{C}$. Then the eigenvalues of $Y_N$ exhibit the following behavior. As $N \to \infty$, for each $1 \le i \le r$,
\begin{equation}\label{E:EigenvalueMappingForward}
\lambda_i (Y_N) \xrightarrow{\text{a.s.}} \eta_i =
\left\{
  \begin{array}{ll}
    \mathcal{G}_{\mu_C}^{-1} (\frac{1}{\theta_i}), & \hbox{if $\theta_i > \frac{1}{ \mathcal{G}_{\mu_C}(a^+)}$} \\
    a, & \hbox{otherwise.}
  \end{array}
\right.
\end{equation}
\end{thm}

\begin{cor}[Eigenvalue phase transition for estimation]\label{T:EigenvalueNonlinearMapping}
Let $Y_N = \sum_{i = 1}^{r} \theta_i \phi_{i} \phi_{i}' + C_N$. Assume that $r \ll N$ and, as $N \to \infty$, $\lambda_{i} (Y_N) \xrightarrow{\text{a.s.}} \eta_i $, $\mu_{C_N} \xrightarrow{\text{a.s.}} \mu_{C}$ (compactly supported). Let $a$ be the supremum of the support of $\mu_{C}$. Then for any $\eta_{i} > a^+$, there exists a corresponding $\theta_i$ such that $\theta_{i} = \frac{1}{\mathcal{G}_{\mu_C} (\eta_{i})}$.
\end{cor}

\begin{rem}
It is worth noting that the supremum of the support of $\mu_{C}$ (denoted by $a$) is actually a critical point for eigenvalue estimation. That means, as $N \to \infty$, for any $1 \le i \le r$, if $\lambda_i (Y_N) < a$, there is indeed no valuable information contained in $\lambda_i (Y_N)$ about the $i$-th eigenvalue of $X_N$ (denoted by $\theta_i$).
\end{rem}

Corollary~\ref{T:EigenvalueNonlinearMapping} yields a method to estimate top eigenvalues associated with a low rank matrix from its noisy measurements. Specifically, for $N$ sufficiently large, if $\lambda_i (Y_N) > a$, $\widehat{\theta_{i}} = \frac{1}{\mathcal{G}_{\mu_C} (\lambda_i (Y_N))}$ is a (nonlinear shrinkage) estimator for the $i$-th largest eigenvalue of $X_N$.

Next, we employ an example in \cite{benaych2011eigenvalues}, where $C_N$ is a Wigner matrix, to illustrate the eigenvalue phase transition phenomena. A Wigner random matrix $V_{N}$ is an Hermitian matrix consisting of i.i.d. random variables with zero mean and unit variance. Benaych-Georges, F. and Nadakuditi, R.R. demonstrated that, for $Y_N = \sum_{i = 1}^{r} \theta_i \phi_{i} \phi_{i}' + \sigma V_N$, as $N \to \infty$,
\begin{equation}\label{E:PhaseTransition}
\lambda_i (Y_N) \xrightarrow{\text{a.s.}} \eta_i =
\left\{
  \begin{array}{ll}
    \theta_i + \frac{\sigma^2}{\theta_i}, & \hbox{if $\theta_i > \sigma$} \\
    2 \sigma, & \hbox{otherwise.}
  \end{array}
\right.
\end{equation}
On the other hand, the proposed estimator for the eigenvalues associated with a low rank matrix from its noisy measurements is
\begin{equation}\label{E:PhaseTransitionEstimation}
\theta_i \xrightarrow{\text{a.s.}}
\left\{
  \begin{array}{ll}
    \frac{\eta_i + \sqrt{\eta_i^2 - 4 \sigma^2}}{2}, & \hbox{if $\eta_i > 2 \sigma$} \\
    \le \sigma, & \hbox{if $\eta_i \le 2 \sigma$.}
  \end{array}
\right.
\end{equation}

\begin{figure}
\centering
\subfigure[Upward bias in the large eigenvalues (eigenvalue phase transition for observation).]{
\label{F:LowRank_WignerShrinkage:1}
\includegraphics[width=0.47\columnwidth]{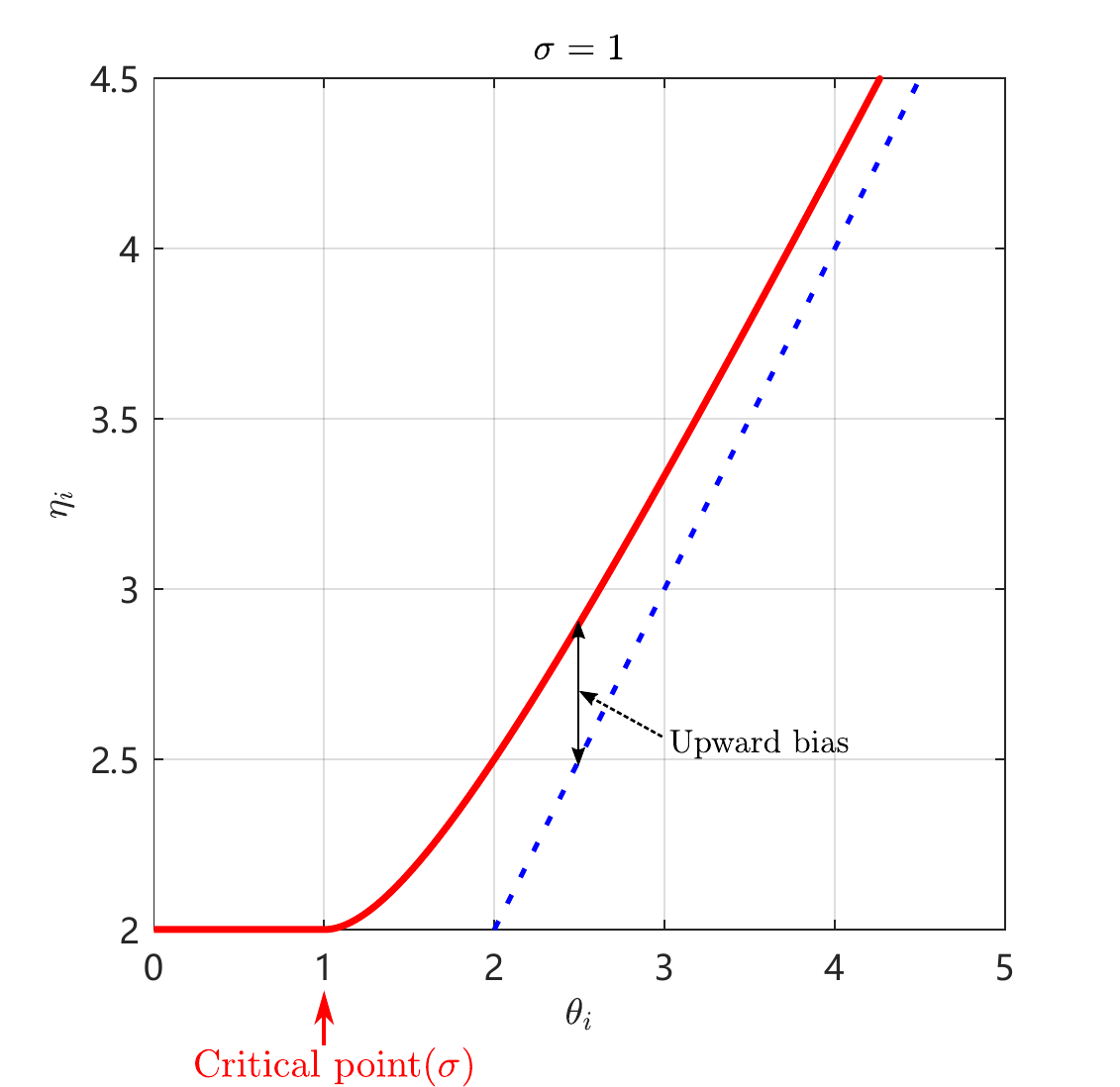}}
\subfigure[Nonlinear shrinkage estimator of the large eigenvalues (eigenvalue phase transition for estimation).]{
\label{F:LowRank_WignerShrinkage:2}
\includegraphics[width=0.47\columnwidth]{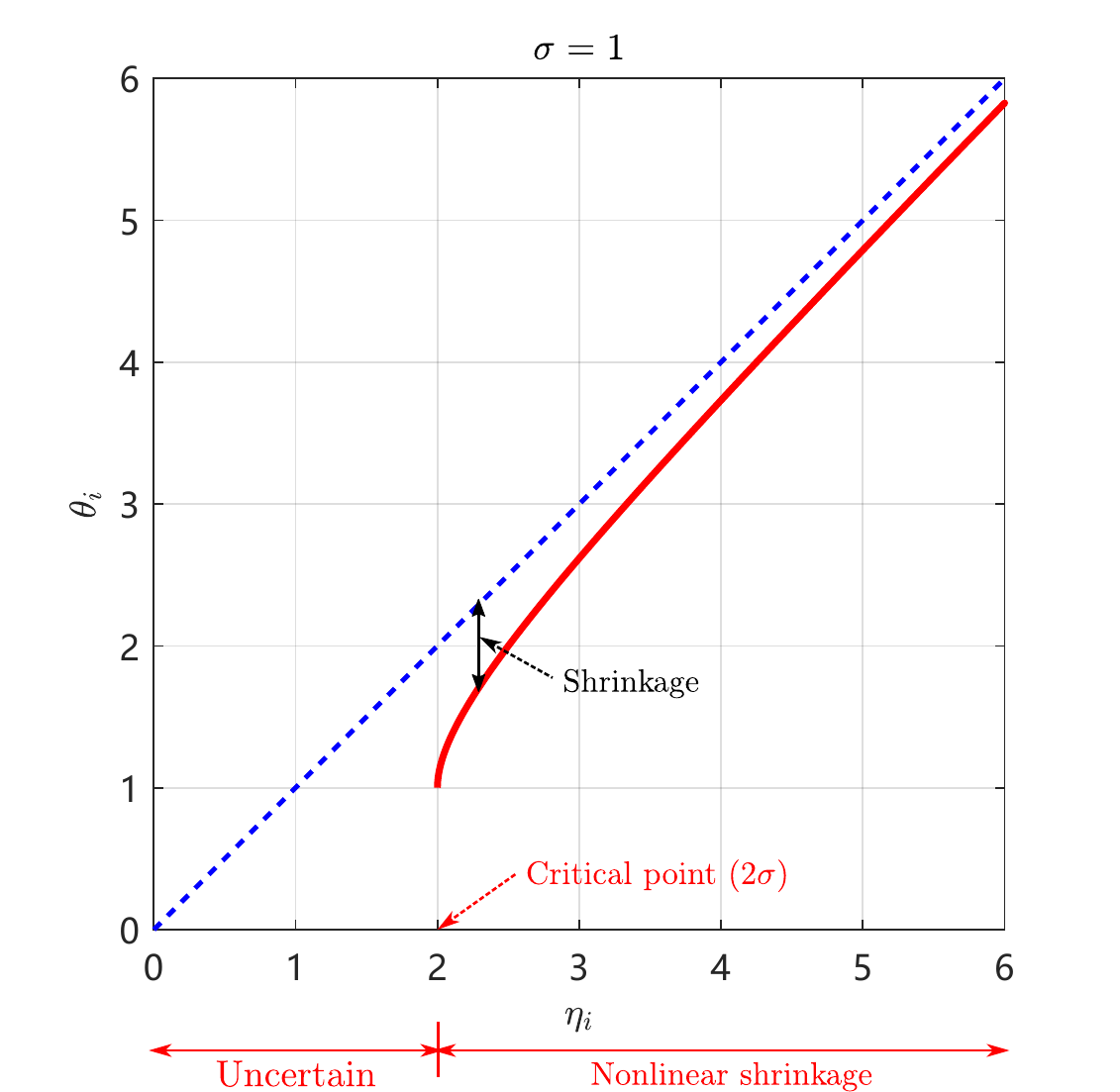}}
\caption{Eigenvalue phase transitions for low-rank spiked sample covariance model. The noise covariance matrix is a Wigner matrix. The mappings defined by  Equation~\ref{E:PhaseTransition} and Equation~\ref{E:PhaseTransitionEstimation} are plotted in red in (a) and (b), respectively. The dashed blue is $y = x$ for easy comparison.}
\label{F:LowRank_WignerShrinkage}
\end{figure}

Fig.~\ref{F:LowRank_WignerShrinkage:1} and \ref{F:LowRank_WignerShrinkage:2} show the mappings defined by Equations~\ref{E:PhaseTransition} and \ref{E:PhaseTransitionEstimation}, respectively. We can see from Fig~\ref{F:LowRank_WignerShrinkage:1} that there exists significant upward bias if the original eigenvalue (low rank) is greater than $\sigma$ (critical point). In Fig.~\ref{F:LowRank_WignerShrinkage:2}, we demonstrate  the nonlinear shrinkage for estimating the large eigenvalues associated with a low rank matrix, where the critical point is $2 \sigma$. That means, for any $1 \le i \le r$, if $\lambda_i (Y_N)$ falls into the uncertain region, i.e., $\lambda_i (Y_N) < 2 \sigma$, there is indeed no chance to recover the $i$-th eigenvalue of $X_N$ (denoted by $\theta_i$) from $\lambda_i (Y_N)$ directly.

\section{Free probability in a nutshell}\label{S:FP}
Free probability is a mathematical theory developed by Voiculescu around 1983. The freeness or free independence is an analogue of the classical notion of independence. A comprehensive study on free probability could be found in \cite{mingo2017free, speicher2014free}.

The significance of free probability to random matrix theory lies in the fundamental observation that random matrices consisting of independent entries also tend to be asymptotically free. Roughly speaking, free probability replaces the vague notion of generic position between eigenspaces associated with random matrices by a mathematical precise concept of freeness and provides a framework to calculate the LSD of $f(A_N, B_N)$ based on LSDs of $A_N$ and $B_N$. Consequently, many tedious calculations in random matrix theory, particularly those of an algebraic or combinatorial nature, can be performed efficiently and systematically with free probability theory \cite{speicher2014free}.

In the following, we will present some preliminaries for free probability theory. For those who are familiar with these may skip most of this section.

We generally refer to a pair $(\mathcal{A}, \phi)$, consisting of a unital algebra $\mathcal{A}$ and a unital linear functional $\phi: \mathcal{A} \to \Cbb$ with $\phi(1) = 1$, as a non-commutative probability space. The elements of $\mathcal{A}$ are always called non-commutative random variables. $\phi(a_{j_1} \cdots a_{j_n})$ for random variables $a_1, \ldots, a_k \in \mathcal{A}$ are always called moments.

\begin{defn}[\cite{mingo2017free}]
Let $(\mathcal{A}, \phi)$ be a unital algebra with a unital linear functional. Suppose $\mathcal{A}_1, \ldots, \mathcal{A}_s$ are unital subalgebras. We say that $\mathcal{A}_1, \ldots, \mathcal{A}_s$ are freely independent (or just free) with respect to $\phi$ if whenever we have $r \ge 2$ and $a_1, \ldots, a_r \in \mathcal{A}$ such that
\begin{itemize}
  \item $\phi(a_i) = 0$ for $i = 1, \ldots, r$
  \item $a_i \in \mathcal{A}_{j_i}$ with $1 \le j_i \le s$ for $i = 1, \ldots, r$
  \item $j_1 \neq j_2, j_2 \neq j_3, \ldots, j_{r-1} \neq j_r$
\end{itemize}
we must have $\phi(a_1 \cdots a_r) = 0$.
\end{defn}

\begin{defn}[\cite{mingo2017free}]\label{D:Freeness}
Let $(\mathcal{A}, \phi)$ be a non-commutative probability space. The non-commutative random variables $a_1, \ldots, a_r \in \mathcal{A}$ are said to be free or freely independent if the generated unital subalgebras $\mathcal{A}_i = \text{alg} (1, a_i)(i=1, \ldots, s)$ are free in $\mathcal{A}$ with respect to $\phi$.
\end{defn}

\begin{rem}
With freeness, it is possible to calculate the mixed moments of $a$ and $b$ from the moments of $a$ and $b$. In this sense, freeness could be considered and treated as the non-commutative analogue to the concept of independence of classical random variables.
\end{rem}

\begin{rem}
Suppose $\mu$ and $\nu$ are two compactly supported probability measures. If there exist two free random variables $a$ and $b$ having distributions $\mu$ and $\nu$, respectively, then the distributions of the non-commutative random variables $a + b$ and $a b$ depend only on $\mu$ and $\nu$ \cite{hiai2006semicircle, couillet2011random}.
\end{rem}

There are many classes of random matrices which show asymptotic freeness. In fact, Gaussian random matrices and deterministic matrices, Haar distributed unitary random matrices and deterministic matrices, Wigner random matrices and deterministic matrices all become asymptotically free with respect to the normalized trace \cite{mingo2017free}. We only present one theorem to demonstrate that randomly rotated deterministic matrices become free in the asymptotic regime.

\begin{thm}[\cite{mingo2017free}]\label{T:HaarRotationFreeness}
Let $A_N$ and $B_N$ be two sequences of deterministic $N \times N$ matrices with $A_N \xrightarrow{distr} a$ and $B_N \xrightarrow{distr} b$. Let $H_N$ be a sequence of $N \times N$ Haar unitary random matrices. Then $A_N$, $H_N B_N H_N^* \xrightarrow{distr} a, b$, where $a$ and $b$ are free. This convergence holds also almost surely. So in particular, we have that $A_N$ and $H_N B_N H_N^*$ are almost surely asymptotically free.
\end{thm}

\subsection{Free harmonic analysis}
We begin with the definition of Cauchy transform associated with a self-adjoint random variable. Actually, the distribution of a self-adjoint random variable $a$ in a $C^*$-probability space $(\mathcal{A}, \phi)$ can be identified by a compactly supported probability measure $\mu_a$ on $\R$, i.e., for all $n \in \mathbb{N}$,
\[
\phi(a^n) = \int_{\R} t^n d \mu_a (t) .
\]
With the Cauchy transform associated with probability measure $\mu_a$, we can then define the Cauchy transform of $a$ as
\[
\mathcal{G}_{a}(z) \equiv \phi \left( \frac{1}{z - a} \right) = \int_{\R} \frac{1}{z-t} d \mu_a (t) .
\]

Next, we will define several transforms on the complex upper half-plane that will be used frequently in this article. The first is the M-transform, defined as
\[
\mathcal{M}_{a}(z) \equiv z \mathcal{G}_{a}(z) - 1.
\]

\begin{lem}\label{T:M_transform}
For a probability measure $\mu$, let $\mathcal{M}_\mu (z)$ be the associated M-transform. Then $\text{Re} \left( \mathcal{M_{\mu}} (i y) \right) < 0$, for all $y > 0$, and $\lim_{y \to \infty}  \mathcal{M_{\mu}} (i y) = 0$. Furthermore, if $ y > 0$ and $\text{supp} (\mu) \in (0, \infty)$, $\text{Im} \left( \mathcal{M}_{\mu} (i y) \right) < 0$.
\end{lem}

\begin{proof}
Since $\mathcal{M}_\mu (i y) = - y  \text{Im} ( \mathcal{G}_\mu (i y) ) - 1 + i y \text{Re} ( \mathcal{G}_\mu (i y) )$, we have that, for $y > 0$,
\[
\begin{aligned}
\text{Re} \left( \mathcal{M}_\mu (i y) \right) &= - y  \text{Im} ( \mathcal{G}_\mu (i y) ) - 1 \\
&= - \int_{\mathbb{R}} y \text{Im} \left( \frac{1}{iy - t} \right) d \mu(t) - 1 \\
&= \int_{\mathbb{R}} \frac{- t^2}{y^2 + t^2} d \mu(t) < 0
\end{aligned}
\]
and
\[
\text{Im} \left( \mathcal{M}_\mu (i y) \right) = y \text{Re} ( \mathcal{G}_\mu (i y) ) = \int_{\mathbb{R}} \frac{-y t}{y^2 + t^2} d \mu(t).
\]
By the assumption that $\text{supp} (\mu) \in (0, \infty)$, $\text{Im} \left( \mathcal{M}_\mu (i y) \right) = \int_{\mathbb{R}^+} \frac{-y t}{y^2 + t^2} d \mu(t) < 0$ for all $y > 0$.
Finally, with the dominated convergence theorem, we have $\lim_{y \to \infty} \text{Re} \left( \mathcal{M}_\mu (i y) \right) = 0$ and $\lim_{y \to \infty} \text{Im} \left( \mathcal{M}_\mu (i y) \right) = 0$, which lead to $\lim_{y \to \infty}  \mathcal{M} (i y) = 0$.
\end{proof}

\begin{thm}[Free multiplication law \cite{burda2011applying, burda2010random}]\label{T:FreeMultiplication}
For a non-commutative random variable $a$, let $\mathcal{M}_{a}(z)$ be the associated M-transform. We define the N-transform $\mathcal{N}_{a} (z)$ as a function such that
$\mathcal{M}_{a} \left( \mathcal{N}_{a} (z) \right) = z$, or equivalently, $\mathcal{N}_{a} \left( \mathcal{M}_{a} (z) \right) = z$. For freely independent random variables $a$ and $b$, the N-transform of $ab$ obeys the `non-commutative multiplication law':
\[
\mathcal{N}_{ab}(z) = \frac{z}{1+z} \mathcal{N}_{a}(z) \mathcal{N}_{b} (z).
\]
\end{thm}

Since the trace operator is invariant under cyclic permutations, there exist similar cyclic properties for the M- and N- transforms.

\begin{thm}[Cyclic property]
Consider an $N \times T$ matrix $A$ and a $T \times N$ matrix $B$ and assume $AB$ is Hermitian. Then
\[
\mathcal{M}_{AB} (z) = \frac{T}{N} \mathcal{M}_{BA} (z),
\]
and
\[
\mathcal{N}_{AB} (z) = \mathcal{N}_{BA} \left( \frac{N}{T} z \right) .
\]
\end{thm}

\begin{proof}
On account of the definition of M-transform, we have $\mathcal{M}_{M} (z) = \sum_{k \ge 1} \frac{\text{tr} \left( M^k \right) }{z^k}$. Since $AB$ and $BA$ have the same non-zero eigenvalues and $AB$ is Hermitian, all of the eigenvalues of $BA$ are real. Then
\[
\begin{aligned}
N \mathcal{M}_{AB} (z) &= N \sum_{k \ge 1} \frac{\text{tr} \left( (AB)^k \right) }{z^k} = \sum_{k \ge 1} \frac{\text{Tr} \left( (AB)^k \right) }{z^k} \\
& = \sum_{k \ge 1} \frac{\text{Tr} \left( (BA)^k \right) }{z^k} = T \sum_{k \ge 1} \frac{\text{tr} \left( (BA)^k \right) }{z^k} \\
& = T \mathcal{M}_{BA} (z) .
\end{aligned}
\]
Applying the scaling property of inverse functions yields
\[
\mathcal{N}_{AB} (z) = \mathcal{N}_{BA} \left( \frac{N}{T} z \right) .
\]
\end{proof}

\section{Separability between the spatio- and temporal- correlations}\label{S:Separability}
As mentioned in Section~\ref{S:Introduction}, the separability brings us significant computational advances. For example, if we address a separable Gaussian process $U = \left( \Sigma^{s}_{N} \right)^{\frac{1}{2}} W_{N,T} \left( \Sigma^{t}_{T} \right)^{\frac{1}{2}}$ as a vector $\text{vec}(U)$ by stacking the columns of $U$ into a single vector, the population covariance matrix $\text{E} \left[ \text{vec}(U) \text{vec}(U)^{\text{T}} \right] = \Sigma^{t}_{T} \otimes \Sigma^{s}_{N}$, where $\otimes$ denotes the Kronecker product between matrices.

\begin{lem}\label{L:SeparableCovarianceMatrix}
Let $U = \left( \Sigma^{s}_{N} \right)^{\frac{1}{2}} W_{N,T} \left( \Sigma^{t}_{T} \right)^{\frac{1}{2}}$, where $W_{N,T}$ is a Gaussian random matrix, $\Sigma^{s}_{N}$ and $\Sigma^{t}_{T}$ are non-negative definite matrices. Then $\text{E} \left[ {\text{vec}(U) \text{vec}(U)^{\text{T}}} \right] = \Sigma^{t}_{T} \otimes \Sigma^{s}_{N}$.
\end{lem}

\begin{proof}
Since $U = \left( \Sigma^{s}_{N} \right)^{\frac{1}{2}} W_{N,T} \left( \Sigma^{t}_{T} \right)^{\frac{1}{2}}$, we first have $\text{vec}(U) = [ ( { \Sigma^{t}_{T} }^{\frac{1}{2}} )^{\text{T}} \otimes { \Sigma^{s}_{N} }^{\frac{1}{2}} ] \text{vec}(W_{N,T})$. Then
\begin{equation*}
\begin{aligned}
& \text{E} \left[ {\text{vec}(U) \text{vec}(U)^{\text{T}}} \right] \\
= & [ ( { \Sigma^{t}_{T} }^{\frac{1}{2}} )^{\text{T}} \otimes { \Sigma^{s}_{N} }^{\frac{1}{2}} ] \text{E} [ \text{vec}(W_{N,T}) \text{vec}(W_{N,T})^{\text{T}} ] [ ( { \Sigma^{t}_{T} }^{\frac{1}{2}} )^{\text{T}} \otimes { \Sigma^{s}_{N} }^{\frac{1}{2}} ]^{\text{T}} \\
= & [ ( { \Sigma^{t}_{T} }^{\frac{1}{2}} )^{\text{T}} \otimes { \Sigma^{s}_{N} }^{\frac{1}{2}} ] [ { \Sigma^{t}_{T} }^{\frac{1}{2}} \otimes ( { \Sigma^{s}_{N} }^{\frac{1}{2}} )^{\text{T}} ] \\
= & \Sigma^{t}_{T} \otimes \Sigma^{s}_{N} .
\end{aligned}
\end{equation*}
\end{proof}

\begin{rem}
The separability between spatio- and temporal-correlations reduces the model complexity greatly. Consider the estimation of covariance matrix as an example, the number of parameters for a separable Gaussian random process is $O(N^2+T^2)$, in contrast to $O(N^2T^2)$ for a general Gaussian random process.
\end{rem}

The objective of this study is to characterize the asymptotic behaviour of the sample covariance matrices $C_N = \frac{1}{T} \left( \Sigma^{s}_{N} \right)^{\frac{1}{2}} W_{N,T} \Sigma^{t}_{T} W_{N,T}^{\text{T}} \left( \Sigma^{s}_{N} \right)^{\frac{1}{2}}$. We place certain restrictions on $\Sigma^{s}_{N}$ and $\Sigma^{t}_{T}$ so that they are completely defined by simple parameters that could be estimated easily. In the remainder of this article, a diagonally dominant Wigner matrix and a shift invariant (Toeplitz) structure are imposed on the spatio-covariance matrix $\Sigma^{s}_{N}$ and the temporal-covariance matrix $\Sigma^{t}_{T}$, respectively.

When $\Sigma^{s}_{N}$ is a diagonally dominant Wigner matrix and $\Sigma^{t}_{T}$ is a deterministic Toeplitz matrix, we could assume the asymptotic freeness between $\frac{1}{T} W_{N,T} \Sigma^{t}_{T} W_{N,T}^{\text{T}}$ and $\Sigma^{s}_{N}$, $\frac{1}{T} W_{N,T}^{\text{T}} W_{N,T}$ and $\Sigma^{t}_{T}$, respectively. The assumption is relatively weak because the left-singular vectors of $W_{N,T}$ compose a Haar-like random matrix. The asymptotic freeness follows immediately from Theorem~\ref{T:HaarRotationFreeness}.

The key idea for the calculation of LSD of $C_N$ is through the non-commutative multiplication law of the N-transform (Theorem~\ref{T:FreeMultiplication}). Specifically, with the asymptotic freeness between $\frac{1}{T} W_{N,T} \Sigma^{t}_{T} W_{N,T}^{\text{T}}$ and $\Sigma^{s}_{N}$, $\frac{1}{T} W_{N,T}^{\text{T}} W_{N,T}$ and $\Sigma^{t}_{T}$, we have
\begin{align*}
& \mathcal{N}_{C_N} (z) \\
= & \mathcal{N}_{ \frac{1}{T} \left( \Sigma^{s}_{N} \right)^{\frac{1}{2}} W_{N,T} \Sigma^{t}_{T} W_{N,T}^{\text{T}} \left( \Sigma^{s}_{N} \right)^{\frac{1}{2}} } (z) \\
= & \mathcal{N}_{ \frac{1}{T}  W_{N,T} \Sigma^{t}_{T} W_{N,T}^{\text{T}} \Sigma^{s}_{N} } (z) \quad \text{(cyclic property)} \\
= & \frac{z}{1+z} \mathcal{N}_{ \frac{1}{T} W_{N,T} \Sigma^{t}_{T} W_{N,T}^{\text{T}} } (z) \mathcal{N}_{ \Sigma^{s}_{N} } (z) \ \text{(free multiplication law)} \\
= & \frac{z}{1+z} \mathcal{N}_{ \frac{1}{T} W_{N,T}^{\text{T}} W_{N,T} \Sigma^{t}_{T} } (\frac{N}{T} z) \mathcal{N}_{ \Sigma^{s}_{N} } (z) \quad \text{(cyclic property)} \\
= & \frac{z}{1+z} \frac{ \frac{N}{T} z }{ 1 + \frac{N}{T} z } \mathcal{N}_{ \frac{1}{T} W_{N,T}^{\text{T}} W_{N,T} } (\frac{N}{T} z) \mathcal{N}_{ \Sigma^{t}_{T} } (\frac{N}{T} z) \mathcal{N}_{ \Sigma^{s}_{N} } (z) \\
& \qquad \qquad \qquad \qquad \qquad \qquad \qquad \ \text{(free multiplication law)} \\
= & \frac{z}{1+z} \frac{ \frac{N}{T} z }{ 1 + \frac{N}{T} z } \mathcal{N}_{ \frac{1}{T} W_{N,T} W_{N,T}^{\text{T}} } ( z ) \mathcal{N}_{ \Sigma^{t}_{T} } (\frac{N}{T} z) \mathcal{N}_{ \Sigma^{s}_{N} } (z). \\
& \qquad \qquad \qquad \qquad \qquad \qquad \qquad \qquad \ \text{(cyclic property)}
\end{align*}

Theorem~\ref{T:MP_Law} shows that, as $N, T \to \infty$ and $\frac{N}{T} \to c \in (0, 1)$, the ESD of $\frac{1}{T} W_{N \times T} W_{N \times T}^{\text{T}}$ converges weakly and almost surely to Mar\v{c}enko-Pastur distribution. By the assumption that $\Sigma^{t}_{T} \xrightarrow{distr} \Sigma^{t}$, $\Sigma^{s}_{N} \xrightarrow{distr} \Sigma^{s}$ and $C_N \xrightarrow{distr} C$, we finally have
\begin{equation}\label{E:N-transformC}
\mathcal{N}_{C} (z) = \frac{z}{1+z} \frac{ c z }{ 1 + c z } \mathcal{N}_{ \text{M-P} } ( z ) \mathcal{N}_{ \Sigma^{t} } (c z) \mathcal{N}_{ \Sigma^{s} } (z) .
\end{equation}

\subsection{Restrictions on the spatio- and temporal- covariance matrices}
Next, we will give a brief exposition of the restrictions on the spatio- and temporal- covariance structures. It is well-known that the autoregressive (AR) model provides flexibility in terms of understanding the lag-dependence structures. Particularly, AR(1) often appears as a building block in more complex models. For an AR(1) Gaussian process, we actually restrict the spatio-covariance matrix $\Sigma^{s}_{N}$ to a diagonally dominant Wigner matrix and impose a `shift invariant' structure to the temporal-covariance matrix $\Sigma^{t}_{T}$, respectively.

Specifically, let $U_{i}(t)$ be a jointly wide-sense stationary AR(1) Gaussian process, i.e., for $1 \le i \le N$,
\[
U_{i} (t) \approx r U_{i} (t-1) + e_{i} (t).
\]
Here $\Abs{r} < 1$, $e_{i}(t)$ is white noise with zero mean. We first calculate the spatio-covariance function
\begin{align*}
\Sigma^{s}_{N} [i, j] \equiv & \text{E} \left[ U_{i}(t) U_{j} (t)  \right] \\
\approx & r^2 \text{E}  \left[ U_{i}(t-1) U_{j}(t-1) \right] + \text{E} \left[ e_{i}(t) e_{j}(t) \right] + \\
&  r \text{E} \left[ U_{i}(t-1) e_{j}(t) \right] + r \text{E} \left[ e_{i}(t) U_{j}(t-1) \right].
\end{align*}
Applying the jointly wide-sense stationarity between $U_i(\cdot)$ and $U_j(\cdot)$ yields
\[
\begin{aligned}
\Sigma^{s}_{N} [i, j] \approx & \frac{\text{E} \left[ e_{i}(t) e_{j}(t) \right]}{1 - r^2} + \\
& \underbrace{\frac{r}{1 - r^2} \left( \text{E} \left[ U_{i}(t-1) e_{j}(t) \right] + \text{E} \left[ e_{i}(t) U_{j}(t-1) \right] \right)}_{\text{cross correlation}}.
\end{aligned}
\]

In practice, there often exists non-trivial cross correlation term $\frac{r}{1 - r^2} ( \text{E} [ U_{i}(t-1) e_{j}(t) ] + \text{E} [ e_{i}(t) U_{j}(t-1) ] )$. Nevertheless, it's not easy to calculate the cross correlation term accurately. To reduce the complexity of the model, we further assume that the random variables in the collection $ \{ \frac{r}{1 - r^2} ( \text{E} [ U_{i}(t-1) e_{j}(t) ] + \text{E} [ e_{i}(t) U_{j}(t-1) ] ) , 1 \le i < j \le N \} $ are i.i.d. with zero mean and standard deviation $\beta$. Additionally, we assume that $\frac{\text{E} \left[ e_{i}(t) e_{j}(t) \right]}{1 - r^2} = \alpha \delta(i - j)$. Finally, the spatio-covariance matrix is modeled as $\Sigma^{s}_{N} = \alpha I_N + \beta V_{N}$, where $V_{N}$ is a Wigner random matrix.

We now turn to the structure of the temporal-covariance matrix $\Sigma^{t}_{T}$ associated with an autoregressive Gaussian process. For a wide-sense stationary process, it is natural to impose `shift invariance' on the temporal-covariance matrix $\Sigma^{t}_{T}$, i.e., the entries of $\Sigma^{t}_{T}$ depend only on the difference between the indices,
\[
\Sigma^{t}_{T} [a, b] = f(a-b), \quad 1 \le a, b \le T.
\]
With the assumption that $\Sigma^{t}_{T} \xrightarrow{\text{distr}} \Sigma^{t}$, $\Sigma^{t}$ becomes a bi-infinite Toeplitz matrix (the indices range over all integers)
\[
\Sigma^{t} [a, b] = f(a-b), \quad a, b \in \mathbb{Z}.
\]

Moreover, for a stationary AR(1) Gaussian process $U_{i} (t) \approx r U_{i} (t-1) + e_{i} (t)$, the covariance between the observations at time $a$ and $b$ is
\begin{align*}
\text{E} \left[ U_{i}(a) U_{i} (b)  \right] & \approx r \text{E} \left[ U_{i}(a-1) U_{i} (b)  \right] + \text{E} \left[ e_{i}(a) U_{i} (b)  \right] \\
& = r \text{E} \left[ U_{i}(a-1) U_{i} (b)  \right].
\end{align*}
It follows immediately that $\Sigma^{t}_{T} [a, b] \approx r \Sigma^{t}_{T} [a-1, b]$. Finally, the bi-infinite Toeplitz matrix becomes a bi-infinite exponential off-diagonal decay matrix, i.e., $\Sigma^{t} [a, b] = r^{ \left| {a-b} \right |}$.

\begin{rem}
We emphasize that the AR process could be considered as a mean-field model to study the behavior of large and complex stochastic models \cite{montanari2018mean, yeo2016random}. That means rather than considering individual signals, the spectrum of the homogeneous covariance matrix is a mean-field approximation of the spectrum of the original heterogenous covariance matrix.
\end{rem}

\begin{figure}
  \centering
  \includegraphics[width=0.4\textwidth]{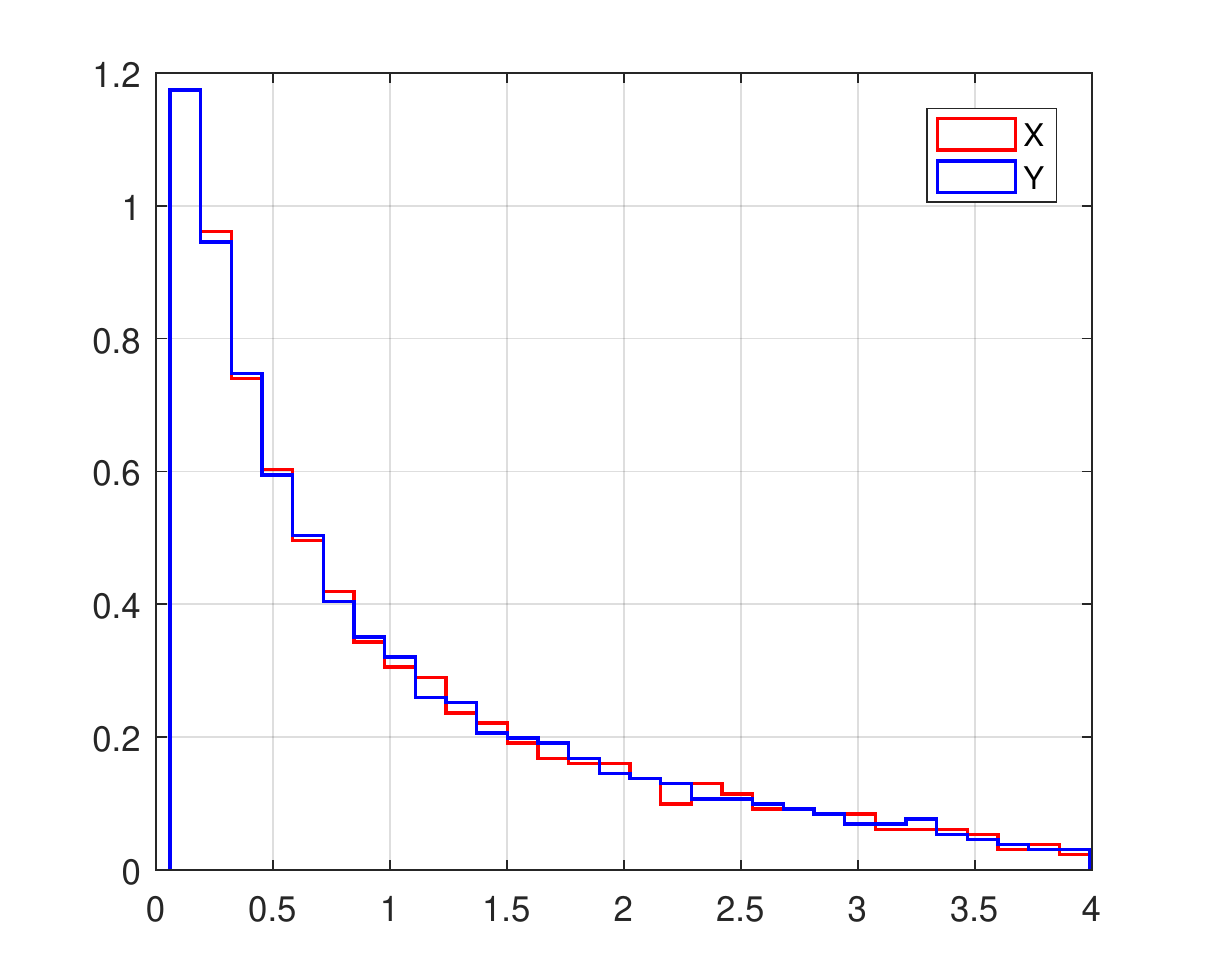}
  \caption{ESDs of the heterogenous (in red) and homogeneous (in blue) autoregressive process.}\label{F:MeanFieldModel}
\end{figure}

We use the mean-field model on spectrum in \cite{yeo2016random} for a better illustration. With a slight abuse of notation, consider two $N \times T$ sample matrices
\[
U[i,t] = r_i U[i,t-1] + \xi[i,t], \quad (\text{heterogenous process})
\]
where $r_i$ is uniformly distributed on $[0, 1]$, $\sigma_i^2 = 1 - r_i^2$ and $\xi[i,t] \sim \mathcal{N}(0, \sigma_i^2)$, and
\[
\bar{U}[i,t] = \bar{r} \bar{U} [i,t-1] + \eta[i,t], \quad (\text{homogeneous process})
\]
where $\bar{r} = \frac{1}{N} \sum r_i$, $\bar{\sigma}^2 = 1 - \bar{r}^2$ and $\eta[i,t] \sim \mathcal{N}(0, \bar{\sigma}^2)$. The ESDs of $\frac{1}{T} U U^T$ and  $\frac{1}{T} \bar{U} \bar{U}^T$ are plotted in Fig.~\ref{F:MeanFieldModel}, from which we can see that there is no significant difference between the two ESDs.

Finally, we have the explicit form of the sample covariance matrix
\[
C_N = \frac{1}{T} \left( \Sigma^{s}_{N} \right)^{\frac{1}{2}} W_{N,T} \Sigma^{t}_{T} W_{N,T}^{\text{T}} \left( \Sigma^{s}_{N} \right)^{\frac{1}{2}},
\]
where $W_{N \times T}$ is a Gaussian random matrix, $\Sigma^{s}_{N} = \alpha I_N + \beta V_{N}$ and $\Sigma^{t}_{T} [a, b] = r^{ \left| {a-b} \right |}$. The aim of this article is to calculate the LSD of $C_N$, under the Kolmogorov condition ($N \to \infty$, $T \to \infty$ and $\frac{N}{T} \to c$).

\section{Cauchy, M- and N-transforms}\label{S:Cauchy_M_N}
We have shown in Section~\ref{S:Separability} that, under the assumption of freeness (relatively weak),
\[
\mathcal{N}_{C} (z) = \frac{z}{1+z} \frac{ c z }{ 1 + c z } \mathcal{N}_{ \text{M-P} } ( z ) \mathcal{N}_{ \Sigma^{t} } (c z) \mathcal{N}_{ \Sigma^{s} } (z).
\]
This section is dedicated to the calculation of Cauchy, M- and N-transforms associated with the Mar\v{c}enko-Pastur distribution, diagonally dominant Wigner matrix ($\Sigma^{s}$) and bi-infinite Toeplitz matrix ($\Sigma^{t}$). It is worth noting that all of these three N-transforms are analytical, so is $\mathcal{N}_{C} (z)$.

Before that, we give a lemma that will be used frequently in the following calculations. For the proof we refer the reader to \cite{mingo2017free}. Unless otherwise stated, we define that $z = \Abs{z} e^{i \theta_z}$, where $0 \le \theta_z < 2 \pi$, and $\sqrt{z} = \sqrt{\Abs{z}} e^{i \frac{\theta_z}{2}}$.

\begin{lem}[\cite{mingo2017free}]\label{Lemma_2}
Let $\zeta_1$ and $\zeta_2$ be the roots of $\zeta^2 - m \zeta + 1 = 0$. Suppose that $\zeta_1 = \frac{m - \sqrt{m^2 - 4}}{2}$ and $\zeta_2 = \frac{m + \sqrt{m^2 - 4}}{2}$. If $m \in \mathbb{C}^+$, then $\zeta_1 \in \text{int}(\Gamma)$ and $\zeta_2 \not \in \text{int}(\Gamma)$. Here $\Gamma= \{\zeta : \Abs{\zeta} = 1 \}$ is the unit circle in complex plane.
\end{lem}

\subsection{Mar\v{c}enko-Pastur distribution}\label{SS:N-transform_MP}
We first calculate the Cauchy, M- and N-transforms associated with the Mar\v{c}enko-Pastur distribution. When $c < 1$, the Mar\v{c}enko-Pastur density function is reduced to
\[
d \mu (x) = \frac{1}{2 \pi c x} \sqrt{(x-a)^{+}(b-x)^{+}}.
\]
Consider the associated Cauchy transform
\begin{equation*}
\mathcal{G}_{\text{M-P}} (z) = \int_{a}^{b} \frac{ \sqrt{(x-a)(b-x)} }{ 2 \pi c x (z -x)} dx.
\end{equation*}
Applying $x = 1 + 2 \sqrt{c} \cos \theta + c \ (0 \le \theta \le \pi)$ yields
\[
\mathcal{G}_{\text{M-P}} (z) = \frac{1}{4 \pi} \int_{0}^{2 \pi} \frac{ 4 \sin^2 \theta }{ \left( 1 + 2 \sqrt{c} \cos \theta + c \right) \left( z - 1 - 2 \sqrt{c} \cos \theta - c \right) } d \theta .
\]
Then make the substitution $\omega = e^{i \theta}$ to write $\mathcal{G}_{\text{M-P}} (z)$ as a contour integral
\[
\mathcal{G}_{\text{M-P}} (z) = \frac{1}{4 c \pi i} \oint_{\Gamma} \frac{ (\omega^2 - 1)^2 }{ \omega \left( \omega^2 + \frac{c + 1}{\sqrt{c}} \omega + 1 \right) \left( \omega^2 - \frac{z - c - 1}{\sqrt{c}} \omega + 1 \right) } d \omega ,
\]
where $\Gamma = \{ \omega : | \omega | = 1 \}$. By Lemma~\ref{Lemma_2} and the residue theorem, we are able to show the closed-form expressions of the Cauchy and M- transforms:
\[
\mathcal{G}_{\text{M-P}} (z) = \frac{z + c -1 - \sqrt{(z - a)(z - b)}}{2 c z} ,
\]
\[
\mathcal{M}_{\text{M-P}} (z) = \frac{z - c - 1 - \sqrt{(z - a)(z - b)}}{2 c}.
\]

Finally, with the inverse relationship $\mathcal{M}_{\text{M-P}} ( \mathcal{N}_{\text{M-P}} (z) ) = z$, we find that $\mathcal{N}_{\text{M-P}} (z)$ satisfies
\[
\frac{\mathcal{N}_{\text{M-P}} (z) - c - 1 - \sqrt{ \left( \mathcal{N}_{\text{M-P}} (z) - a \right) \left( \mathcal{N}_{\text{M-P}} (z) - b \right) }}{2c} = z
\]
Solving this equation will give rise to
\[
\mathcal{N}_{\text{M-P}} (z) = \frac{(1 + z) (1 + c z)}{z} .
\]

\subsection{Diagonally dominant Wigner matrices}
We are now in a position to calculate the Cauchy, M- and N-transforms associated with the shifted semi-circle distribution. The well-known semi-circle law is presented first for completeness.

\begin{thm}[Semi-circle law \cite{bai2010spectral}]\label{T:SemiCircleLaw}
Let $V_{N}$ be an $N \times N$ Wigner matrix. Then, with probability 1, the ESD of $V_{N}$ converges to a non-random distribution function $\mu_V$ given by:
\[
d \mu_V (x) = \left\{
           \begin{array}{ll}
             \frac{1}{2 \pi} \sqrt{4 - x^2}, & \hbox{if $|x| \le 2$;} \\
             0, & \hbox{otherwise.}
           \end{array}
         \right.
\]
\end{thm}

\begin{prop}[Shifted semi-circle law]\label{T:ShiftSemiCircleLaw}
Let $\Sigma^{s}_{N} = \alpha I_N + \beta V_{N}$, where $V_{N}$ is a Wigner matrix, $\alpha \ge 0$ and $\beta > 0$. Then, with probability 1, the ESD of $\Sigma^{s}_{N}$ converges to a non-random distribution function $\mu_{\Sigma^{s}}$ given by:
\[
d \mu_{\Sigma^{s}} (x) = \left\{
           \begin{array}{ll}
             \frac{1}{2 \pi \beta^2} \sqrt{4 \beta^2 - (x-\alpha)^2}, & \hbox{if $\left| x - \alpha \right| \le  2 \beta$;} \\
             0, & \hbox{otherwise.}
           \end{array}
         \right.
\]
\end{prop}

\begin{proof}
On account of the spectral mapping theorem, we have that, for any $\alpha \ge 0$ and $\beta > 0$, $\lambda_{i} \left( \alpha I_N + \beta V_{N} \right) = \alpha + \beta \lambda_{i} \left( V_{N} \right)$. The shifted semi-circle law then follows from the semi-circle law.
\end{proof}

\begin{figure}
\centering
\subfigure[Shifted semi-circle distribution (in red) and ESD (in blue).]{
\label{F:ShiftSemiCircleLaw:1}
\includegraphics[width=0.472\columnwidth]{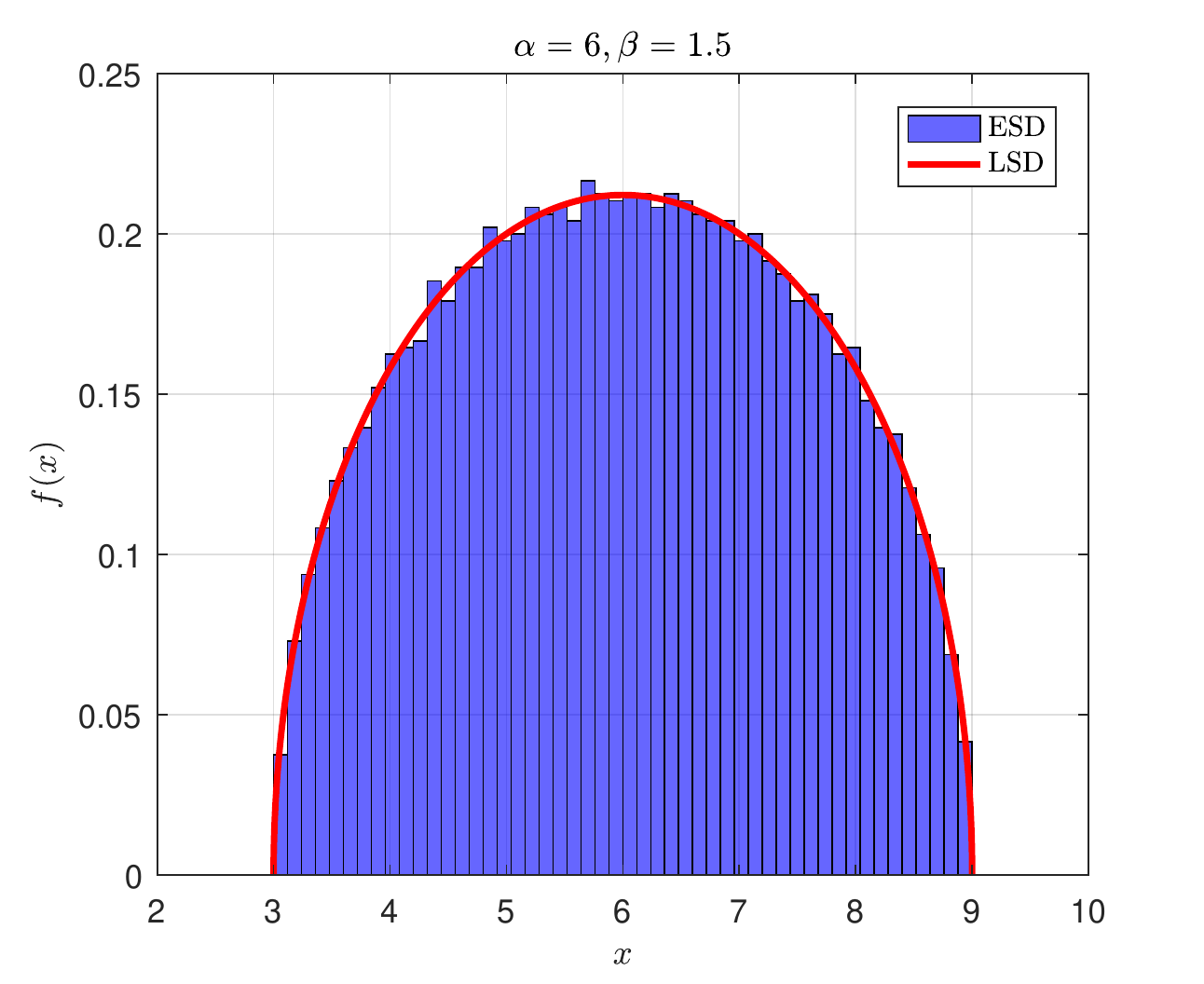}}
\subfigure[Shifted semi-circle distribution for different $\alpha$ and $\beta$.]{
\label{F:ShiftSemiCircleLaw:2}
\includegraphics[width=0.47\columnwidth]{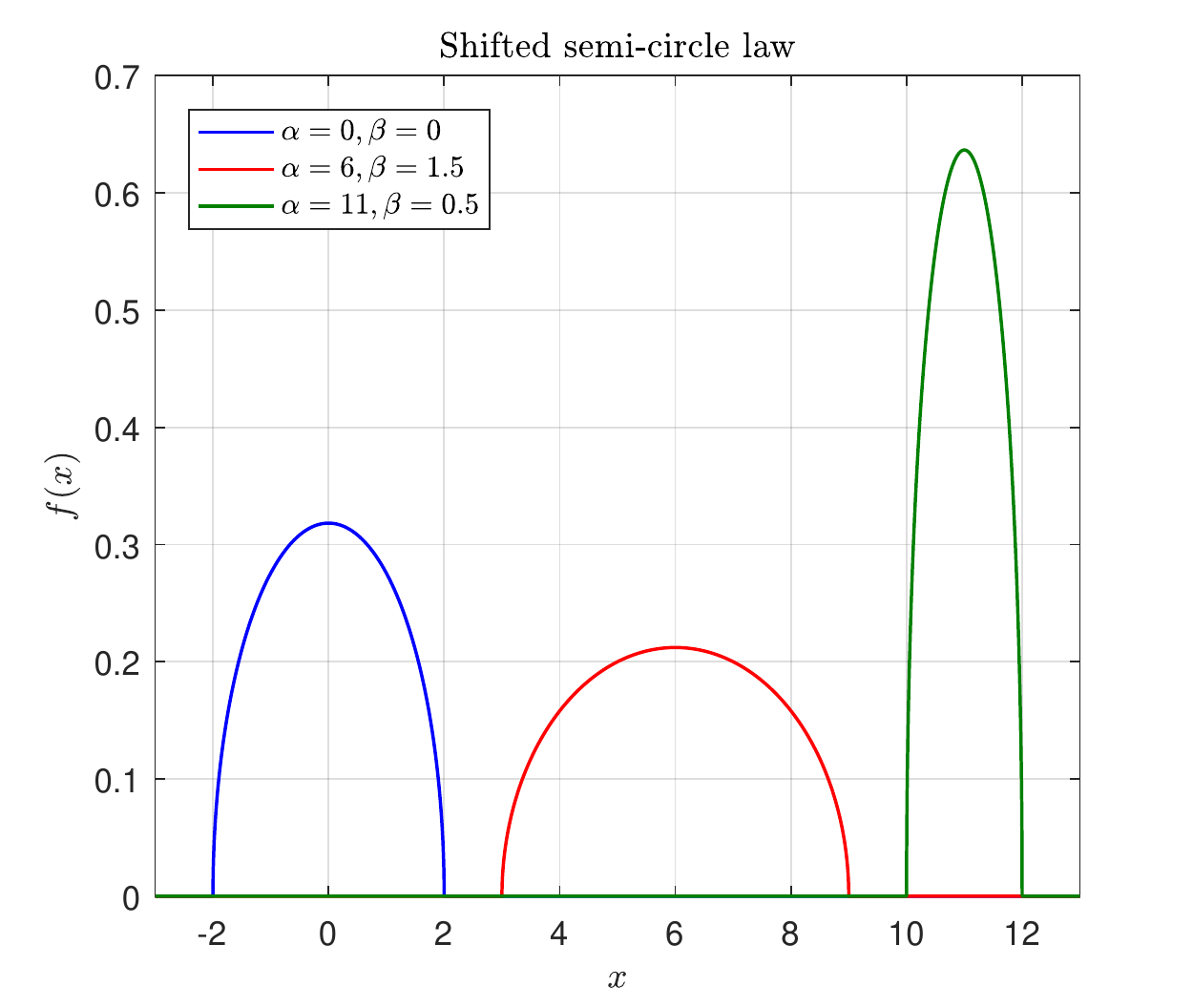}}
\caption{Shifted semi-circle law for different $\alpha$ and $\beta$.}
\label{F:ShiftSemiCircleLaw}
\end{figure}

Proposition~\ref{T:ShiftSemiCircleLaw} asserts that the shifted semi-circle distribution is supported on $\left[ \alpha - 2 \beta, \alpha + 2 \beta \right]$. Fig.~\ref{F:ShiftSemiCircleLaw} shows an excellent agreement between the shifted semi-circle distribution (in red) and the ESD (in blue) of a diagonally dominant Wigner matrix with N = 4000 (one realization).

The Cauchy transform associated with the shifted semi-circle distribution can be calculated in much the same manner as the Mar\v{c}enko-Pastur distribution. Specifically, it follows immediately from Proposition~\ref{T:ShiftSemiCircleLaw} that
\begin{align*}
\mathcal{G}_{\Sigma^s}(z)
& = \frac{1}{2 \pi \beta^2} \int_{\alpha - 2 \beta}^{\alpha + 2 \beta} \frac{\sqrt{4 \beta^2 - (t-\alpha)^2}}{z - t} dt \\
& = \frac{1}{2 \pi \beta^2} \int_{-2\beta}^{2\beta} \frac{\sqrt{4 \beta^2 - t^2}}{z - t - \alpha} dt .
\end{align*}
Applying $t = 2 \beta \cos \theta$ ($0 \le \theta \le \pi$) yields
\[
\mathcal{G}_{\Sigma^s}(z) = \frac{1}{4 \pi} \int_{0}^{2 \pi} \frac{4 \sin^2 \theta}{z - 2 \beta \cos \theta - \alpha} d \theta.
\]
Then make the substitution $\omega = e^{i \theta}$ to write $\mathcal{G}_{\Sigma^s}$ as a contour integral
\[
\mathcal{G}_{\Sigma^s}(z) = \frac{1}{4 \pi \beta i} \oint_{\Gamma} \frac{ (\omega^2 - 1)^2 }{ \omega^2 \left( \omega^2 - \frac{z - \alpha}{\beta} \omega + 1 \right) } d \omega ,
\]
where $\Gamma = \{ \omega : | \omega | = 1 \}$. By Lemma~\ref{Lemma_2} and the residue theorem, we are able to show the closed-form expressions of the Cauchy and M- transforms
\[
\mathcal{G}_{\Sigma^s}(z) = \frac{z - \alpha - \sqrt{(z - \alpha)^2 - 4 \beta^2} }{2 \beta^2},
\]
\[
\mathcal{M}_{\Sigma^s}(z) = \frac{z \left( z - \alpha - \sqrt{(z - \alpha)^2 - 4 \beta^2} \right) }{2 \beta^2} - 1 .
\]

With the inverse relationship $\mathcal{M}_{\Sigma^s}(\mathcal{N}_{\Sigma^s} (z)) = z$, we find that $\mathcal{N}_{\Sigma^s} (z)$ satisfies a quadratic equation
\[
z \mathcal{N}_{\Sigma^s} (z)^2 - \alpha(z + 1) \mathcal{N}_{\Sigma^s} (z) - \beta^2 (z + 1)^2 = 0 .
\]
Solving this equation will lead to
\[
\mathcal{N}_{\Sigma^s} (z) = \frac{(z + 1) \left( \alpha \pm \sqrt{\alpha^2 + 4 \beta^2 z } \right) }{2 z} .
\]
On account of Lemma~\ref{T:M_transform}, we see that, for all $y > 0$, $\text{Re} \left( \mathcal{M} (i y) \right) < 0$ and $\lim_{y \to \infty} \mathcal{M} (i y) = 0 $. Therefore,
\[
\begin{aligned}
\Abs{ \mathcal{N}_{\Sigma^s} ( 0 ) } = & \Abs{ \mathcal{N}_{ \Sigma^s} \left( \lim_{y \to \infty} \mathcal{M}_{\Sigma^s} (i y) \right) } \\
= & \lim_{y \to \infty} \Abs{ \mathcal{N}_{ \Sigma^s} \left(  \mathcal{M}_{\Sigma^s} (i y) \right) } = + \infty .
\end{aligned}
\]
Finally, we have the closed-form expression of the N-transform
\[
\mathcal{N}_{\Sigma^s} (z) = \frac{ (z + 1) \left( \alpha + \sqrt{\alpha^2 + 4 \beta^2 z} \right) }{2 z}.
\]

\subsection{Infinite shift invariant (Toeplitz) matrices}
We now turn to the calculation of Cauchy, M- and N-transforms associated with a bi-infinite Toeplitz matrix. A standard technique to deal with infinite matrices is via the two-dimensional discrete-time Fourier transform. Specifically, we can perform Fourier transform on the row and column indices to produce a function with frequency variables $\omega_1, \omega_2 \in \left[ 0, 2 \pi \right)$,
\[
\hat{A}(\omega_1, \omega_2) = \sum_{a \in \mathbb{Z}} \sum_{b \in \mathbb{Z}} e^{- i (a \omega_1 + b \omega_2)} A [a, b]
\]
or conversely, for $a, b \in \mathbb{Z}$,
\[
A[a, b] = \frac{1}{4 \pi^2} \int_{0}^{2 \pi} \int_{0}^{2 \pi} e^{ i (a \omega_1 + b \omega_2)} \hat{A}( \omega_1, \omega_2) d \omega_1 d \omega_2.
\]

For a bi-infinite Toeplitz matrix $\Sigma^{t} [a, b] = f(a-b)$ and an infinite matrix $G$, the discrete-time Fourier transform of the product $G \Sigma^{t}$ is
\begin{align*}
& \widehat{G \Sigma^{t}} (\omega_1, \omega_2) \\
= & \sum_{a \in \mathbb{Z}} \sum_{b \in \mathbb{Z}} e^{- i (a \omega_1 + b \omega_2)} \sum_{\tau \in \mathbb{Z}} G[a, \tau] \Sigma^{t}[\tau, b] \\
= & \sum_{\tau \in \mathbb{Z}} \sum_{a \in \mathbb{Z}} \sum_{b \in \mathbb{Z}} e^{- i (a \omega_1 + b \omega_2)} G[a, \tau] f(\tau - b) \\
= & \sum_{\tau \in \mathbb{Z}} \sum_{a \in \mathbb{Z}} e^{- i a \omega_1} G[a, \tau] \sum_{b \in \mathbb{Z}} e^{- i b \omega_2} f(\tau - b) \\
= & \sum_{\tau \in \mathbb{Z}} \sum_{a \in \mathbb{Z}} e^{- i (a \omega_1 + \tau \omega_2)} G[a, \tau] \sum_{b \in \mathbb{Z}} e^{ i (\tau - b) \omega_2} f(\tau - b) \\
= & \hat{G} (\omega_1, \omega_2) \hat{f} ( - \omega_2) .
\end{align*}

Next, we will consider the matrix $G = (z I - \Sigma^{t})^{-1}$. It is easy to check that $z \hat{G} (\omega_1, \omega_2) - \widehat{G \Sigma^{t}} (\omega_1, \omega_2) = \hat{I} (\omega_1, \omega_2)$. Substituting $\widehat{G \Sigma^{t}} (\omega_1, \omega_2)$ with $\hat{G} (\omega_1, \omega_2) \hat{f} ( - \omega_2)$  yields
\begin{equation*}
\hat{G} (\omega_1, \omega_2) = \frac{\hat{I} (\omega_1, \omega_2)}{z - \hat{f} (- \omega_2)}.
\end{equation*}
A trivial verification shows that the Fourier transform maps an infinity identity matrix $I[a, b] = \delta (a - b)$ to a Dirac delta function, i.e., $\hat{I}(\omega_1, \omega_2) = 2 \pi \delta(\omega_1 + \omega_2)$. We then have that
\[
\hat{G} (\omega_1, \omega_2) = \frac{2 \pi \delta(\omega_1 + \omega_2)}{z - \hat{f} (- \omega_2)}.
\]
It follows immediately that
\begin{equation}\label{E:G_Matrix_Toeplitz}
\begin{aligned}
G[a, b] = & \frac{1}{4 \pi^2} \int_{0}^{2 \pi} \int_{0}^{2 \pi} d \omega_1 d \omega_2 e^{ i (a \omega_1 + b \omega_2)} \hat{G} (\omega_1, \omega_2) \\
= & \frac{1}{2 \pi} \int_{0}^{2 \pi} \int_{0}^{2 \pi} d \omega_1 d \omega_2 e^{ i (a \omega_1 + b \omega_2)} \frac{\delta(\omega_1 + \omega_2)}{z - \hat{f} (- \omega_2)} \\
= & \frac{1}{2 \pi} \int_{0}^{2 \pi} d \omega e^{ i \omega (a - b)} \frac{1}{z - \hat{f}(- \omega)} .
\end{aligned}
\end{equation}

We can see easily from Equation~\ref{E:G_Matrix_Toeplitz} that $G = (z I - \Sigma^{t})^{-1}$ is actually a shift invariant (Toeplitz) matrix. Applying the normalized trace to the infinite matrix $G$ yields that
\[
\mathcal{G}_{\Sigma^{t}} (z) = \text{tr} (G) = G[0, 0] = \frac{1}{2 \pi} \int_{0}^{2 \pi} \frac{1}{z - \hat{f}(- \omega)} d \omega ,
\]
and
\begin{equation}\label{E:M_Toeplitz}
\mathcal{M}_{\Sigma^{t}} (z) = \frac{1}{2 \pi} \int_{0}^{2 \pi}  \frac{ \hat{f} ( - \omega) }{z - \hat{f}( - \omega) } d \omega .
\end{equation}

We are now in a position to calculate the Cauchy, M- and N-transforms associated with the bi-infinite exponential off-diagonal decay matrix $\Sigma^{t} [a, b] = f(a-b) = r^{ \left| {a-b} \right |}$. It is easy to calculate that
\[
\hat{f} ( - \omega) = \sum_{t \in \mathbb{Z}} e^{i t \omega} f(t) = \frac{1}{1 - r e^{i \omega}} + \frac{1}{1 - r e^{- i \omega}} - 1 .
\]
Substituting $\hat{f} ( - \omega)$ into \ref{E:M_Toeplitz} yields
\begin{align*}
\mathcal{M}_{\Sigma^{t}} (z) = & \frac{1}{2 \pi} \int_{0}^{2 \pi} \frac{ \hat{f} (\omega) }{z - \hat{f}(\omega)} d \omega \\
= & \frac{1}{2 \pi} \int_{0}^{2 \pi} \frac{\frac{1}{1 - r e^{i \omega}} + \frac{1}{1 - r e^{- i \omega}} - 1}{z - (\frac{1}{1 - r e^{i \omega}} + \frac{1}{1 - r e^{- i \omega}} - 1)} d \omega \\
= & \frac{1}{2 \pi} \int_{0}^{2 \pi} \frac{1 - r^2}{ z ( r^2 + 1 ) + r^2 - 1 - z r ( e^{-i \omega} + e^{i \omega} )  } d \omega .
\end{align*}
Make the substitution $\zeta = e^{i \omega}$ to write $\mathcal{M}_{\Sigma^{t}} (z)$ as a contour integral
\begin{align*}
\mathcal{M}_{\Sigma^{t}} (z) = & \frac{1}{2 \pi i} \oint_{\Gamma} \frac{1 - r^2}{- z r \zeta^2 + (z ( r^2 + 1 ) + r^2 - 1) \zeta - z r } d \zeta \\
= & \frac{r^2 - 1}{2 \pi z r i} \oint_{\Gamma} \frac{1}{\zeta^2 - \frac{z ( r^2 + 1 ) + r^2 - 1}{z r} \zeta + 1 } d \zeta,
\end{align*}
where $\Gamma= \{\zeta : \Abs{\zeta} = 1 \}$.

An easy computation shows that, for all $z \in \mathbb{C}^+$ and $0 < r < 1$, $\text{Im} \left( \frac{z ( r^2 + 1 ) + r^2 - 1}{z r} \right) > 0$. By Lemma~\ref{Lemma_2} and the residue theorem, we are able to show the closed-form expressions of M- and Cauchy transforms
\[
\mathcal{M}_{\Sigma^{t}} (z) = \frac{1}{\sqrt{ z - \frac{1+r}{1-r} } \sqrt{ z - \frac{1-r}{1+r} } },
\]
and
\begin{equation}\label{E:G_Topelitz}
\mathcal{G}_{\Sigma^{t}} (z) = \frac{1}{z \sqrt{ z - \frac{1+r}{1-r} } \sqrt{ z - \frac{1-r}{1+r} } } + \frac{1}{z}.
\end{equation}
We finally give a theorem to describe the LSD associated with an exponential off-diagonal decay matrix via the Stieltjes inversion formula.

\begin{thm}\label{T:LSD_Toeplitz}
Let $\Sigma^{t}_T$ be a $T \times T$ exponential off-diagonal decay matrix, i.e., $\Sigma^{t}_T [a, b] = r^{ \left| {a-b} \right |}$, $0 < r < 1$. As $T \to \infty$, the spectral distribution of $\Sigma^{t}_T$ converges to a distribution function $\mu_{\Sigma^{t}}$ given by:
\[
d \mu_{\Sigma^{t}} (x) = \left\{
           \begin{array}{ll}
             \frac{1}{\pi x \sqrt{ (x - a) (b - x)} }, & \hbox{if $x \in (a, b)$;} \\
             0, & \hbox{otherwise.}
           \end{array}
         \right.
\]
where $a = \frac{1-r}{1+r}$ and $b = \frac{1+r}{1-r}$.
\end{thm}

\begin{figure}
\centering
\subfigure[LSD (in red) and ESD (in blue), where $T = 4000$ and $r = 0.5$.]{
\label{F:LSD_ESD_AR1:1}
\includegraphics[width=0.47\columnwidth]{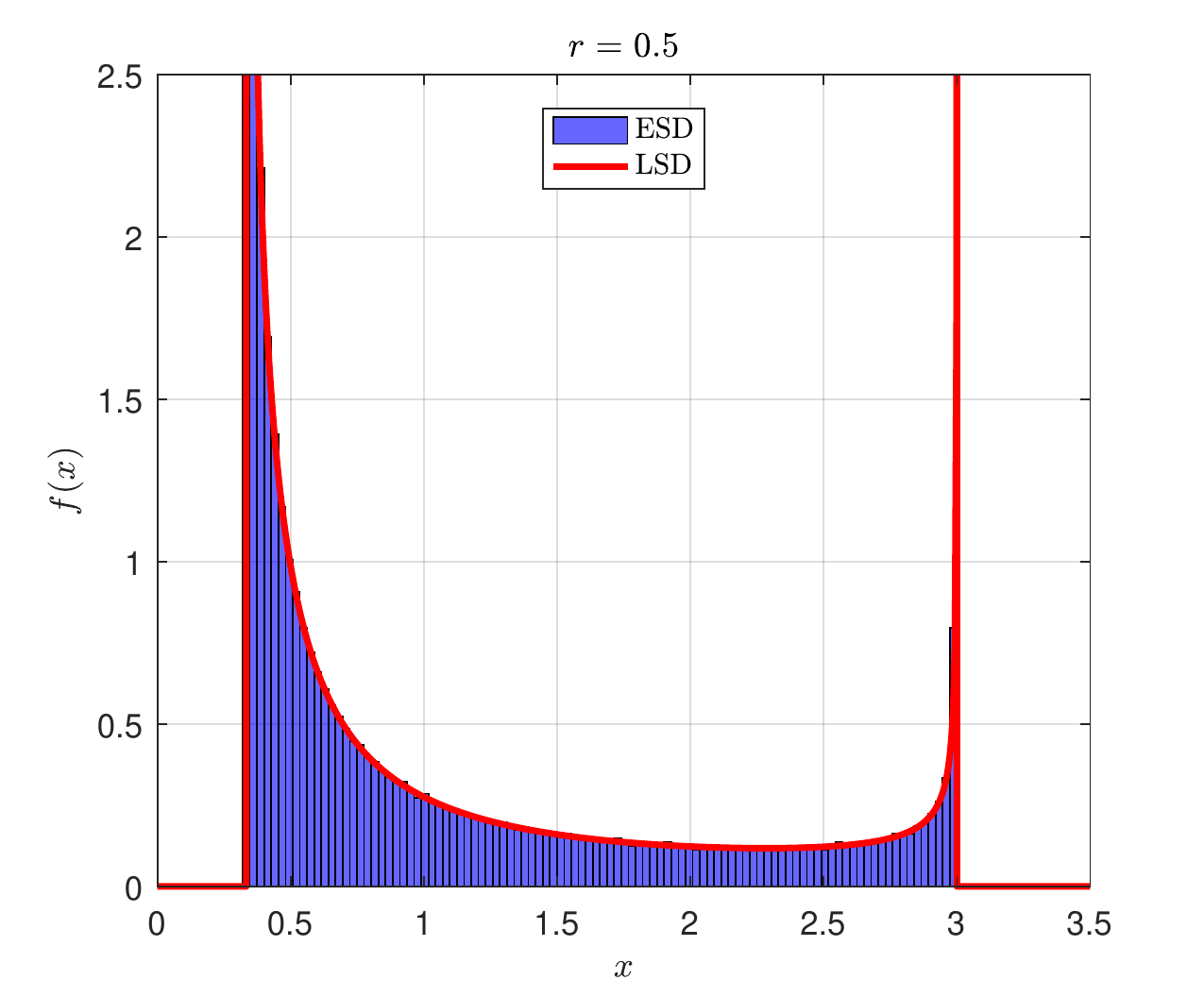}}
\subfigure[LSDs for different values of r.]{
\label{F:LSD_ESD_AR1:2}
\includegraphics[width=0.47\columnwidth]{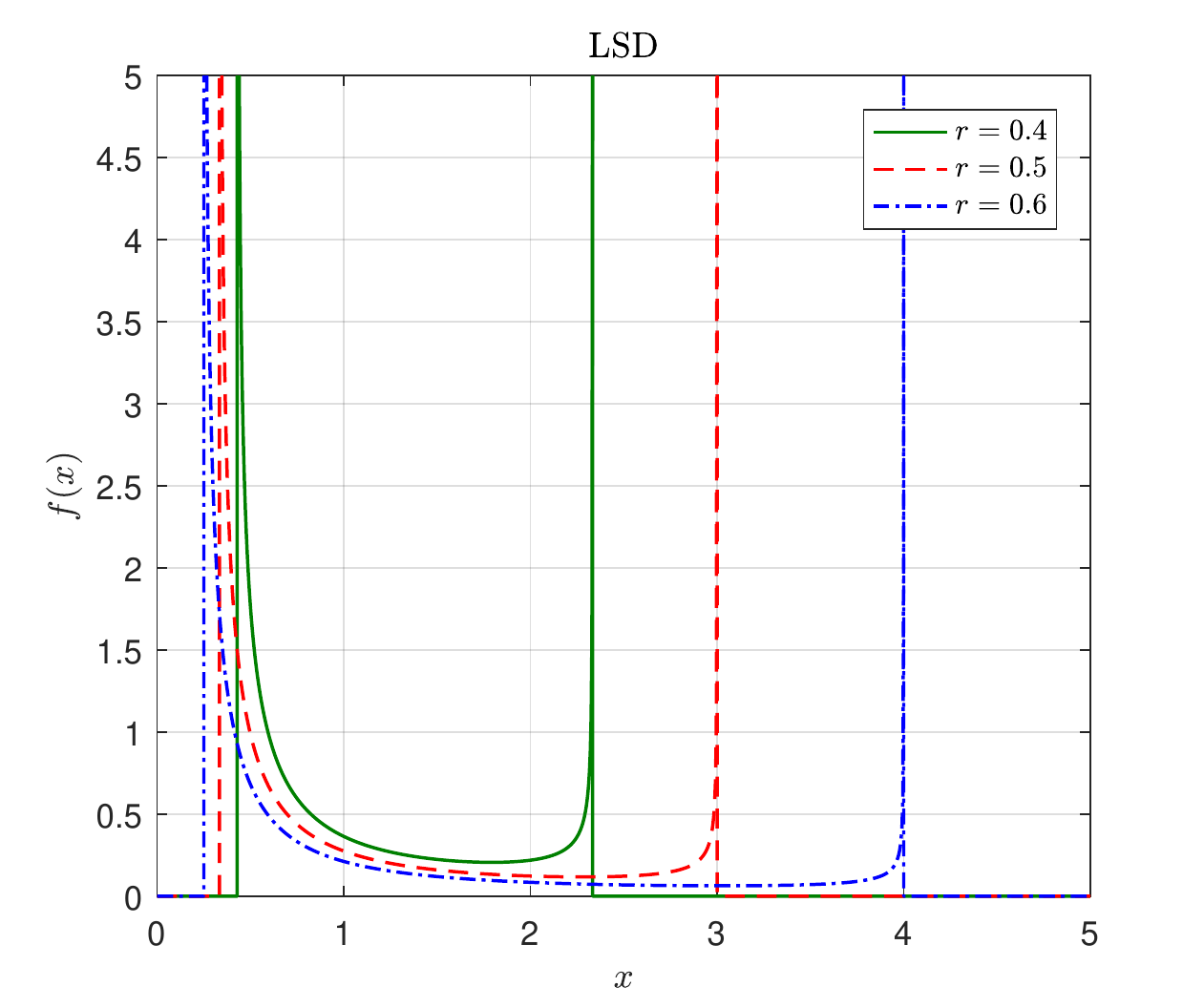}}
\caption{Limiting spectrum associated with exponential off-diagonal decay matrix.}
\label{F:LSD_ESD_AR1}
\end{figure}

For a better illustration of Theorem~\ref{T:LSD_Toeplitz}, we show in Fig.~\ref{F:LSD_ESD_AR1:1} the theoretical LSD (in red) and the spectral distribution of the matrix $\Sigma^{t}_{T}$ (in blue) with $T = 4000$ and $r = 0.5$ (one realization). There is an excellent agreement between the LSD and the spectral distribution. Fig.~\ref{F:LSD_ESD_AR1:2} plots the LSDs for different values of r. We can see that the support of the theoretical LSD spreads out as the value of $r$ increases.

Finally, with the inverse function $\mathcal{M}_{\Sigma^{t}} (\mathcal{N}_{\Sigma^{t}} (z)) = z$, we find that $\mathcal{N}_{\Sigma^{t}} (z)$ satisfies
\[
\mathcal{N}_{\Sigma^{t}}(z)^2 - 2 \frac{1+r^2}{1-r^2} \mathcal{N}_{\Sigma^{t}}(z) + (1 - \frac{1}{z^2}) = 0 .
\]
Solving this yields
\[
\mathcal{N}_{\Sigma^t} (z) = \frac{1+r^2}{1-r^2} \pm \sqrt{\left( \frac{1+r^2}{1-r^2} \right)^2 + \frac{1}{z^2} - 1}.
\]
Furthermore, due to $\mathcal{M}_{\Sigma^{t}} (z) \in \mathbb{C}^-$ for any $z \in \mathbb{C}^+$, we have
\[
\mathcal{N}_{\Sigma^t} (z) = \frac{1+r^2}{1-r^2} + \sqrt{\left( \frac{1+r^2}{1-r^2} \right)^2 + \frac{1}{z^2} - 1}.
\]

\section{Spectrum of large sample covariance matrices of autoregressive Gaussian processes}\label{S:LSD_AR}
The task now is to calculate the LSD associated with a sample covariance matrix
\[
C_N = \frac{1}{T} \left( \Sigma^{s}_{N} \right)^{\frac{1}{2}} W_{N,T} \Sigma^{t}_{T} W_{N,T}^{\text{T}} \left( \Sigma^{s}_{N} \right)^{\frac{1}{2}},
\]
where $W_{N \times T}$ is a Gaussian random matrix, $\Sigma^{s}_{N} = \alpha I_N + \beta V_{N}$ and $\Sigma^{t}_{T} [a, b] = r^{ \left| {a-b} \right |}$. Instead of obtaining the associated M- and Cauchy transforms directly, our approach actually relies on finding the canonical equations that the transforms satisfy.

For the convenience of the reader, we first list the N-transforms associated with the Mar\v{c}enko-Pastur distribution, a diagonally dominant Wigner matrix $\Sigma^{s}_{N} = \alpha I_N + \beta V_{N}$ and a bi-infinite exponential off-diagonal decay matrix $\Sigma^{t}_{T} [a, b] = r^{ \left| {a-b} \right |}$
\[
\mathcal{N}_{\text{M-P}} (z) = \frac{(1 + z) (1 + c z)}{z}
\]
\[
\mathcal{N}_{\Sigma^s} (z) = \frac{ (z + 1) \left( \alpha + \sqrt{\alpha^2 + 4 \beta^2 z} \right) }{2 z}
\]
\[
\mathcal{N}_{\Sigma^t} (c z) = \gamma + \sqrt{\gamma^2 + \frac{1}{c^2 z^2} - 1}
\]
where $\gamma = \frac{1+r^2}{1-r^2}$. Substituting $\mathcal{N}_{\text{M-P}} (z)$, $\mathcal{N}_{\Sigma^s} (z)$ and $\mathcal{N}_{\Sigma^t} (c z)$ into
\[
\mathcal{N}_{C} (z) = \frac{z}{1+z} \frac{ c z }{ 1 + c z } \mathcal{N}_{ \text{M-P} } ( z ) \mathcal{N}_{ \Sigma^{t} } (c z) \mathcal{N}_{ \Sigma^{s} } (z)
\]
yields that
\[
\begin{aligned}
& \mathcal{N}_{C} (z) \\
= & \frac{c}{2} ( z + 1 ) \left( \gamma + \sqrt{\gamma^2 - 1 + \frac{1}{c^2 z^2}} \right) \left( \alpha + \sqrt{\alpha^2 + 4 \beta^2 z} \right) .
\end{aligned}
\]
With the inverse function $\mathcal{M}_{\Sigma^{t}} (\mathcal{N}_{\Sigma^{t}} (z)) = z$ and $\mathcal{M}_{a}(z) = z \mathcal{G}_{a}(z) - 1$, we find that $\mathcal{M}_{C} (z)$ and $\mathcal{G}_{C}(z)$ satisfy
\[
\begin{aligned}
z = & \frac{c}{2} \left( \mathcal{M}_{C} (z) + 1 \right) \left( \gamma + \sqrt{\gamma^2 - 1 + \frac{1}{c^2 \mathcal{M}_{C} (z) ^2}} \right) \\
& \left( \alpha + \sqrt{\alpha^2 + 4 \beta^2 \mathcal{M}_{C} (z)} \right)
\end{aligned}
\]
and
\begin{equation}\label{E:G_Order8a}
\begin{aligned}
1 = & \frac{c}{2} \mathcal{G}_{C}(z) \left( \gamma + \sqrt{ \gamma^2 - 1 + \frac{1}{c^2 (z \mathcal{G}_{C}(z) - 1)^2} }  \right) \\
& \left( \alpha + \sqrt{\alpha^2 + 4 \beta^2 \left( z \mathcal{G}_{C}(z) - 1 \right) } \right)
\end{aligned}
\end{equation}
respectively. Furthermore, with some algebras, Equation~\ref{E:G_Order8a} could be written as an 8-th order polynomial
\begin{equation}\label{E:G_Order8b}
\sum_{i = 0}^{8} p_i \mathcal{G}_{C} (z) ^i = 0,
\end{equation}
where
\begin{equation*}
\begin{aligned}
p_8 = & \beta^4 z^4 , \\
p_7 = & -4 \beta^4 z^3 , \\
p_6 = & \beta^4 z^2 - \frac{2 \beta^4 z^2}{c^2} + \frac{2 \alpha \beta^2 \gamma z^3}{c} , \\
p_5 = & -4 \beta^4 z + \frac{4 \beta^4 z}{c^2} - \frac{6 \alpha \beta^2 \gamma z^2}{c} + \frac{2 \beta^2 z^3}{c^2} - \frac{4 \beta^2 \gamma^2 z^3}{c^2} , \\
p_4 = & \beta^4 + \frac{\beta^4}{c^4} - \frac{2 \beta^4}{c^2} - \frac{2 \alpha \beta^2 \gamma z}{c^3} + \frac{6 \alpha \beta^2 \gamma z}{c} + \frac{\alpha^2 z^2}{c^2} \\
& - \frac{6 \beta^2 z^2}{c^2} + \frac{12 \beta^2 \gamma^2 z^2}{c^2} , \\
p_3 = & \frac{2 \alpha \beta^2 \gamma}{c^3} - \frac{2 \alpha \beta^2 \gamma}{c} - \frac{2 \beta^2 z}{c^4} - \frac{2 \alpha^2 z}{c^2} + \frac{6 \beta^2 z}{c^2} \\
& - \frac{12 \beta^2 \gamma^2 z}{c^2} - \frac{2 \alpha \gamma z^2}{c^3} , \\
p_2 = & -\frac{\alpha^2}{c^4} + \frac{2 \beta^2}{c^4} + \frac{\alpha^2}{c^2} - \frac{2 \beta^2}{c^2} + \frac{4 \beta^2 \gamma^2}{c^2} + \frac{4 \alpha \gamma z}{c^3} + \frac{z^2}{c^4} , \\
p_1 = & -\frac{2 \alpha \gamma}{c^3} - \frac{2 z}{c^4} , \\
p_0 = & \frac{1}{c^4} .
\end{aligned}
\end{equation*}

In much the same way as the calculations in Section~\ref{S:Cauchy_M_N}, in order to calculate the LSD, we need to determine the roots of the polynomial in $\mathcal{G}_{C}$ firstly and then use the Stieltjes inversion formula. Nevertheless, unlike the situations in Section~\ref{S:Cauchy_M_N}, it seems hopeless to find the closed-form expression of $\mathcal{G}_{C}$ by solving Equation~\ref{E:G_Order8b} directly. We finally give a semi-closed-form expression of the LSD. Specifically, since the support of the LSD is unknown, we could conjecture such a region and solve numerically Equation~\ref{E:G_Order8b} for every $z$. The Stieltjes inversion formula can then be used to obtain the LSD.

\begin{figure}
\centering
\subfigure[Theoretical LSD (in red) and ESD (in blue), $N = 2000$, $T = 4000$.]{
\label{F:LSD_ESD_C:1}
\includegraphics[width=0.47\columnwidth]{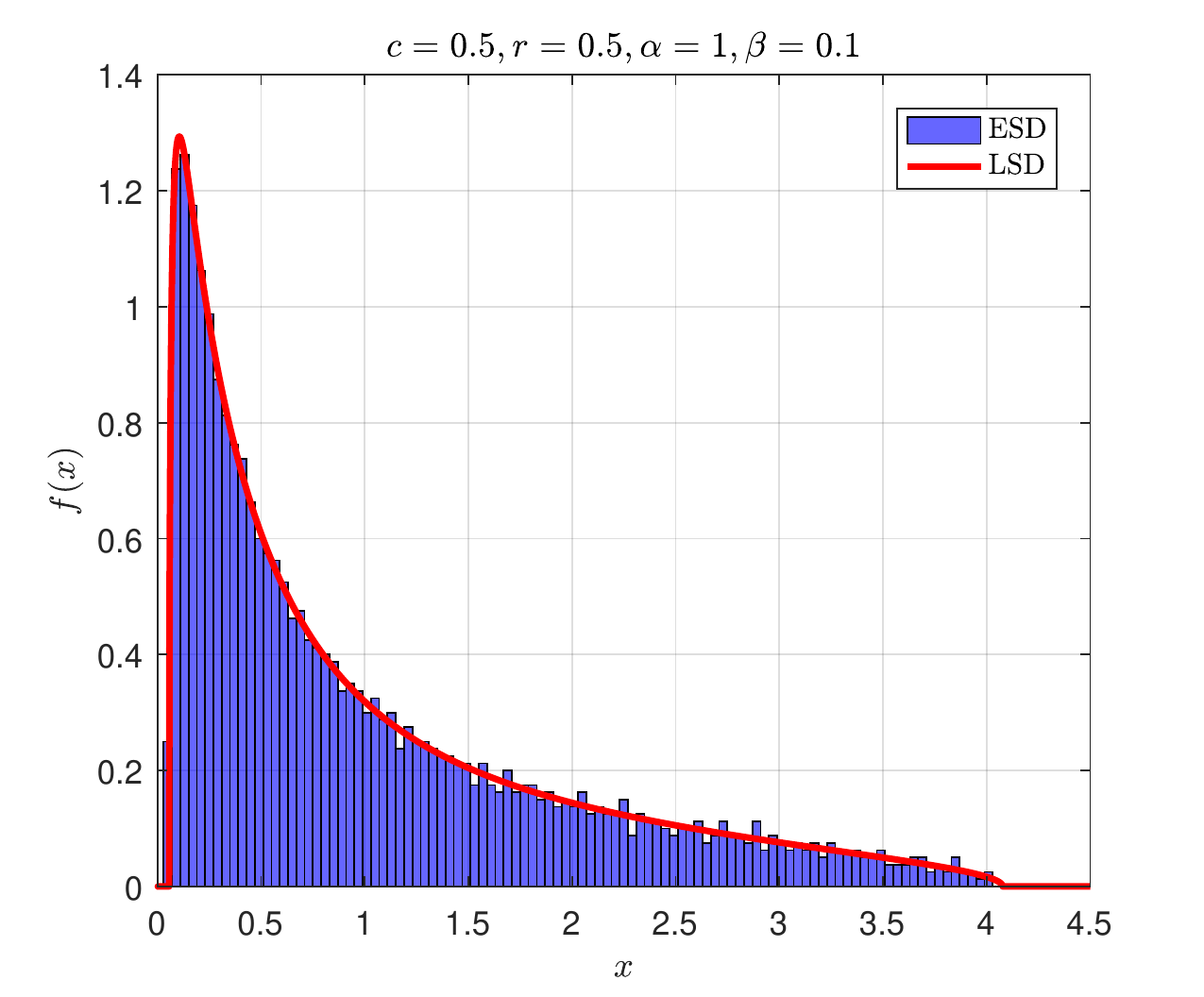}}
\subfigure[LSDs for different values of $r$.]{
\label{F:LSD_ESD_C:2}
\includegraphics[width=0.47\columnwidth]{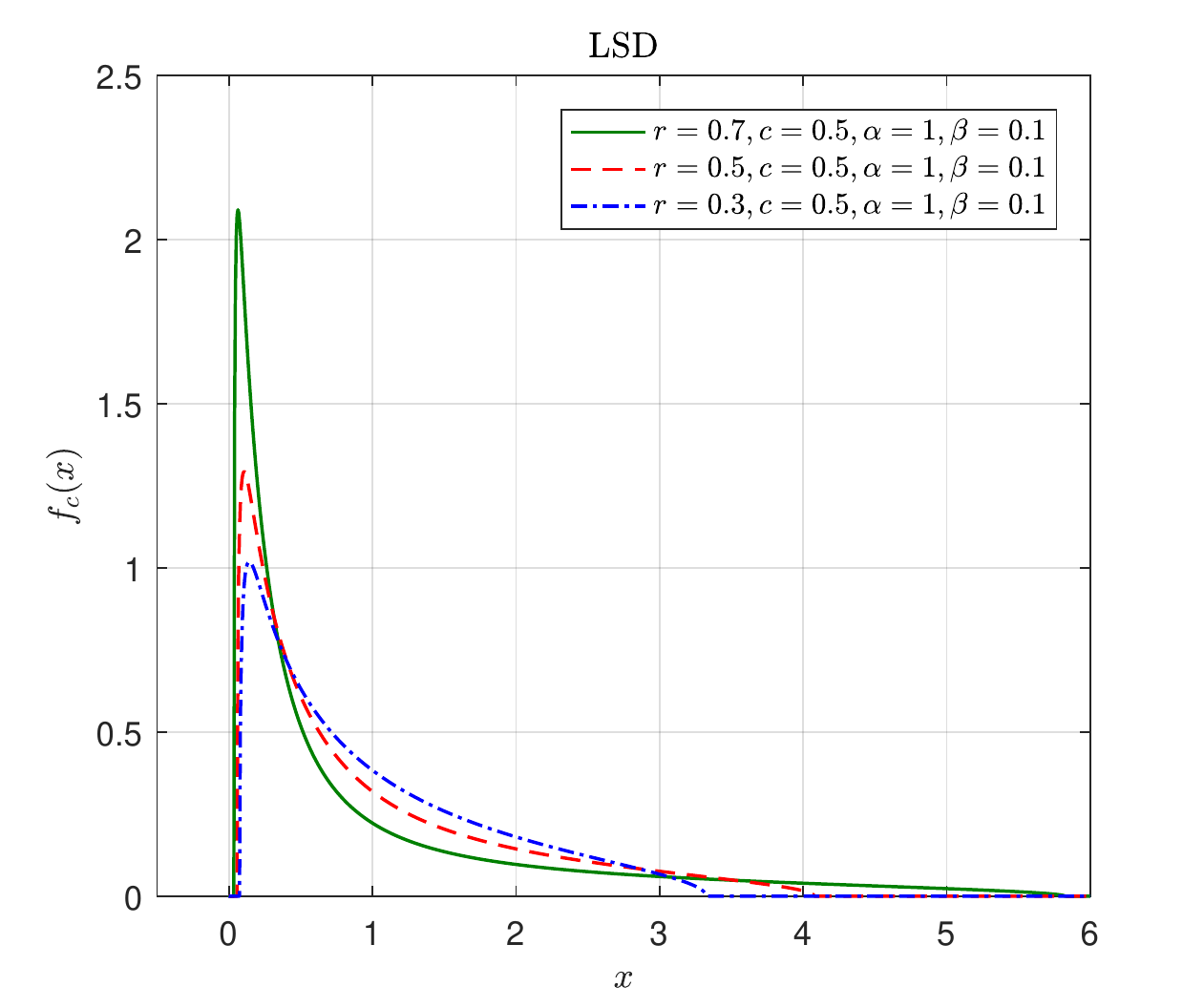}}
\newline
\subfigure[LSDs for different values of $c$.]{
\label{F:LSD_ESD_C:3}
\includegraphics[width=0.47\columnwidth]{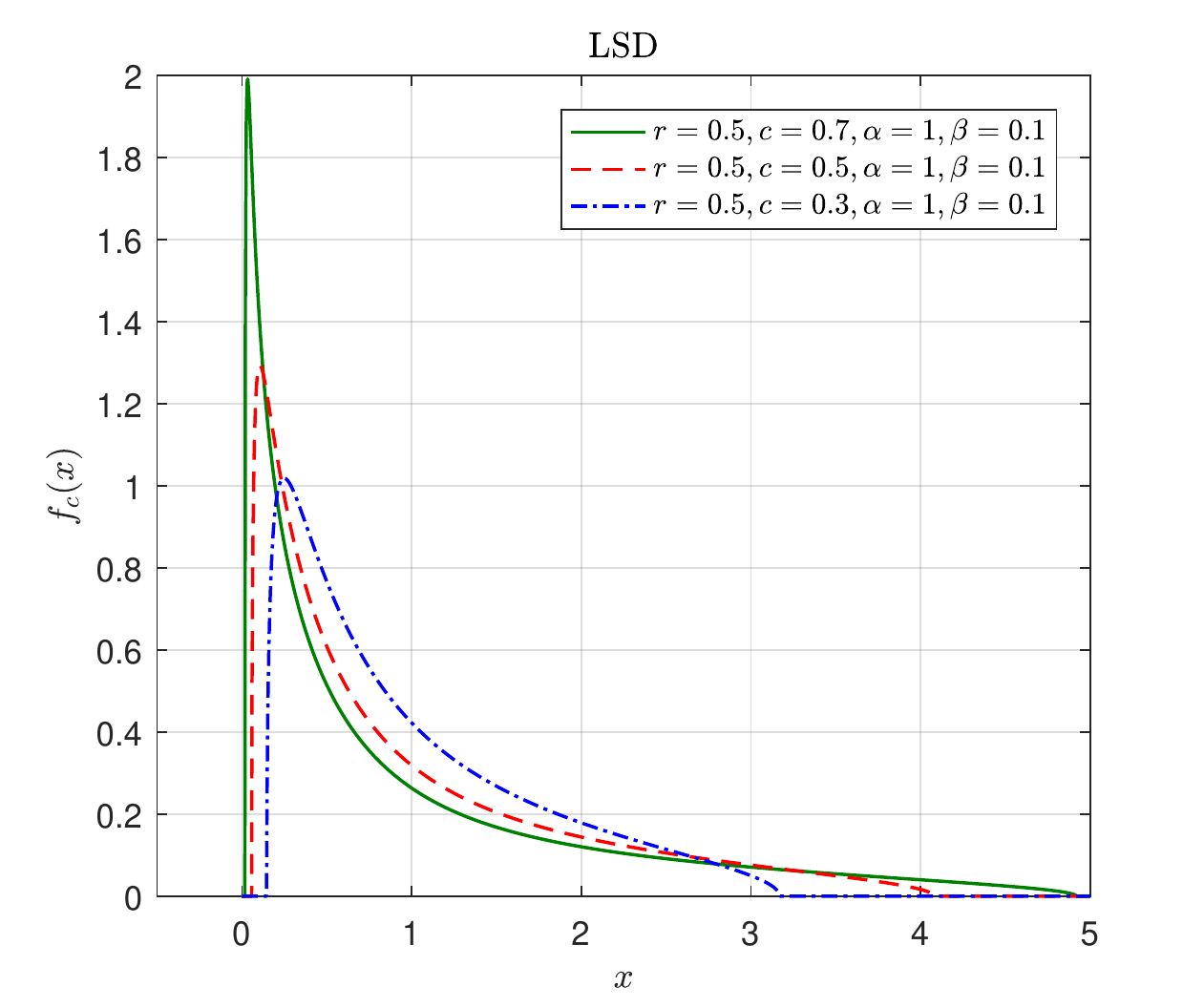}}
\subfigure[LSDs for different values of $\beta$.]{
\label{F:LSD_ESD_C:4}
\includegraphics[width=0.47\columnwidth]{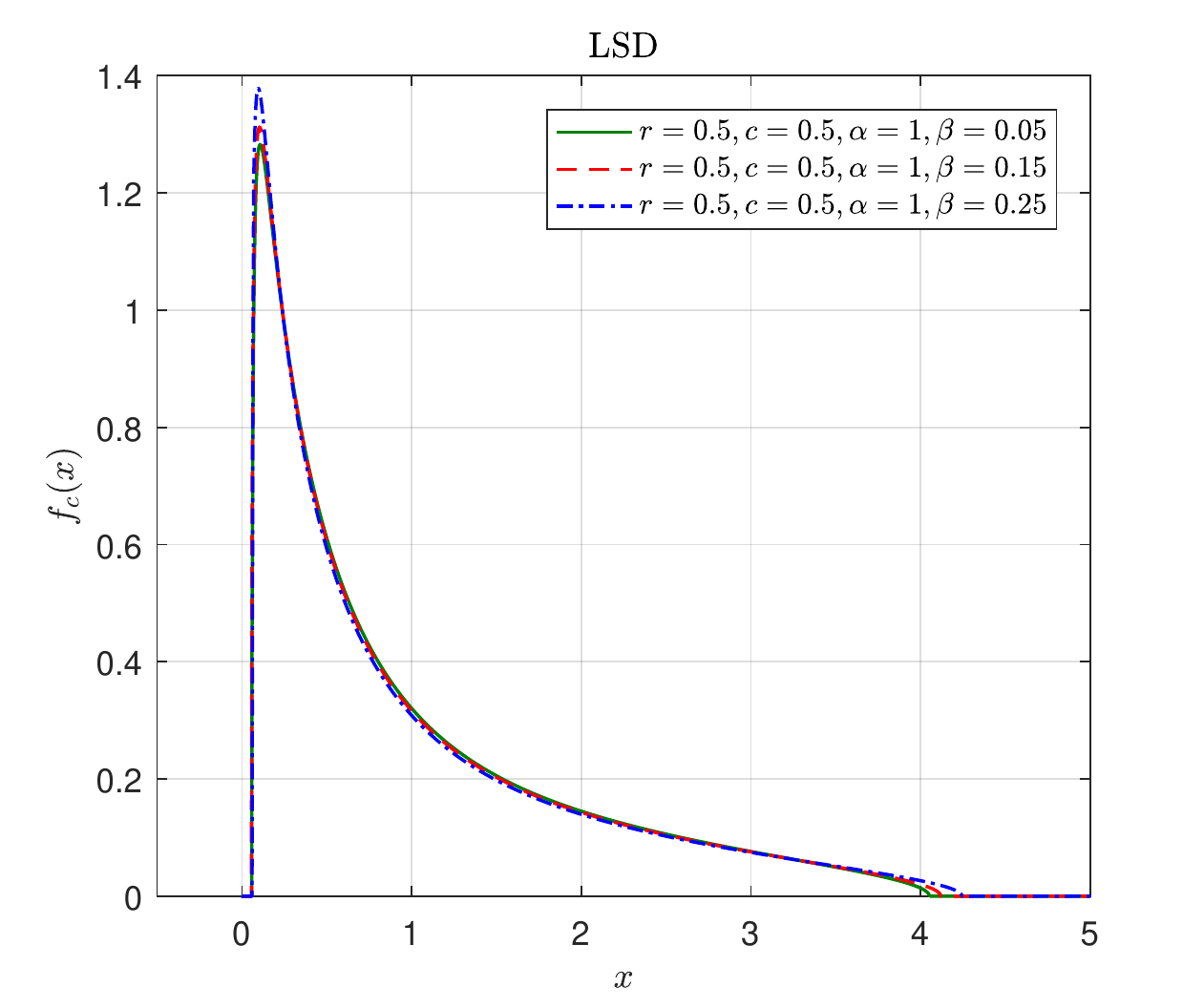}}
\caption{Theoretical LSD and ESD of $C_N = \frac{1}{T} \left( \Sigma^{s}_{N} \right)^{\frac{1}{2}} W_{N,T} \Sigma^{t}_{T} W_{N,T}^{\text{T}} \left( \Sigma^{s}_{N} \right)^{\frac{1}{2}}$.}
\label{F:LSD_ESD_C}
\end{figure}

Fig.~\ref{F:LSD_ESD_C:1} shows an excellent agreement between the theoretical LSD (in red) obtained from numerically solving Equation~\ref{E:G_Order8b} and the ESD (in blue) of $C_N$ with $N = 2000$, $T = 4000$ (one realization). Comparing the distribution in Fig.~\ref{F:LSD_ESD_C:1} with the Mar\v{c}enko-Pastur distribution in Fig.~\ref{F:MP_Law:1}, we can see that the support of the distribution in Fig.~\ref{F:LSD_ESD_C:1} becomes much larger. Fig.~\ref{F:LSD_ESD_C:2}, Fig.~\ref{F:LSD_ESD_C:3} and Fig.~\ref{F:LSD_ESD_C:4} show the behaviors of the LSDs for different values of $r$, $c$ and $\beta$, while other parameters are fixed. We may notice from Fig.~\ref{F:LSD_ESD_C:2} and Fig.~\ref{F:LSD_ESD_C:3} that the support of the LSD spreads out as the value of $c$ or $r$ decreases.

\begin{figure}
\centering
\subfigure[Heterogenous process]{
\label{F:LSD_ESD_AR_Residual:1}
\includegraphics[width=0.47\columnwidth]{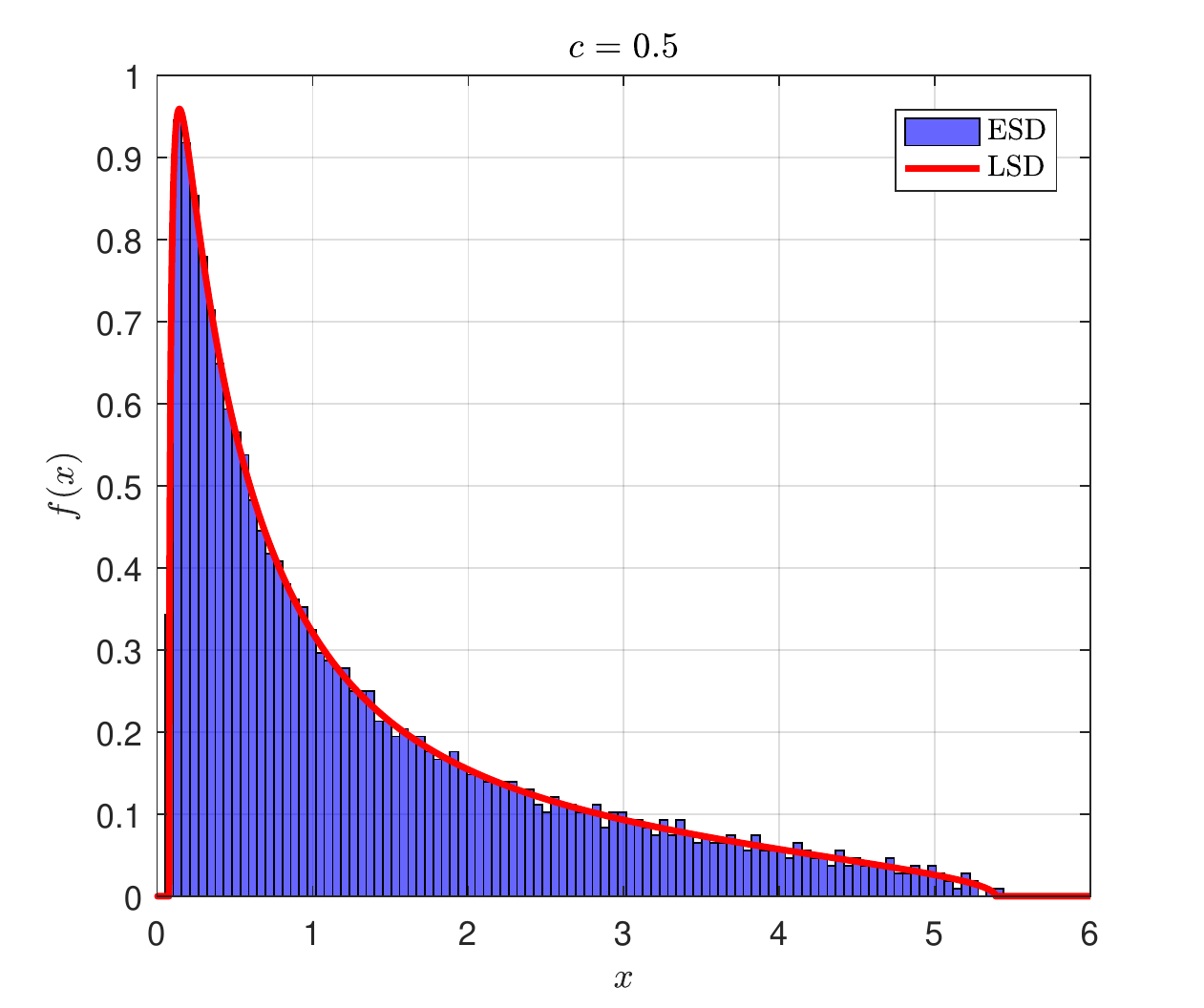}}
\subfigure[Homogeneous process]{
\label{F:LSD_ESD_AR_Residual:2}
\includegraphics[width=0.47\columnwidth]{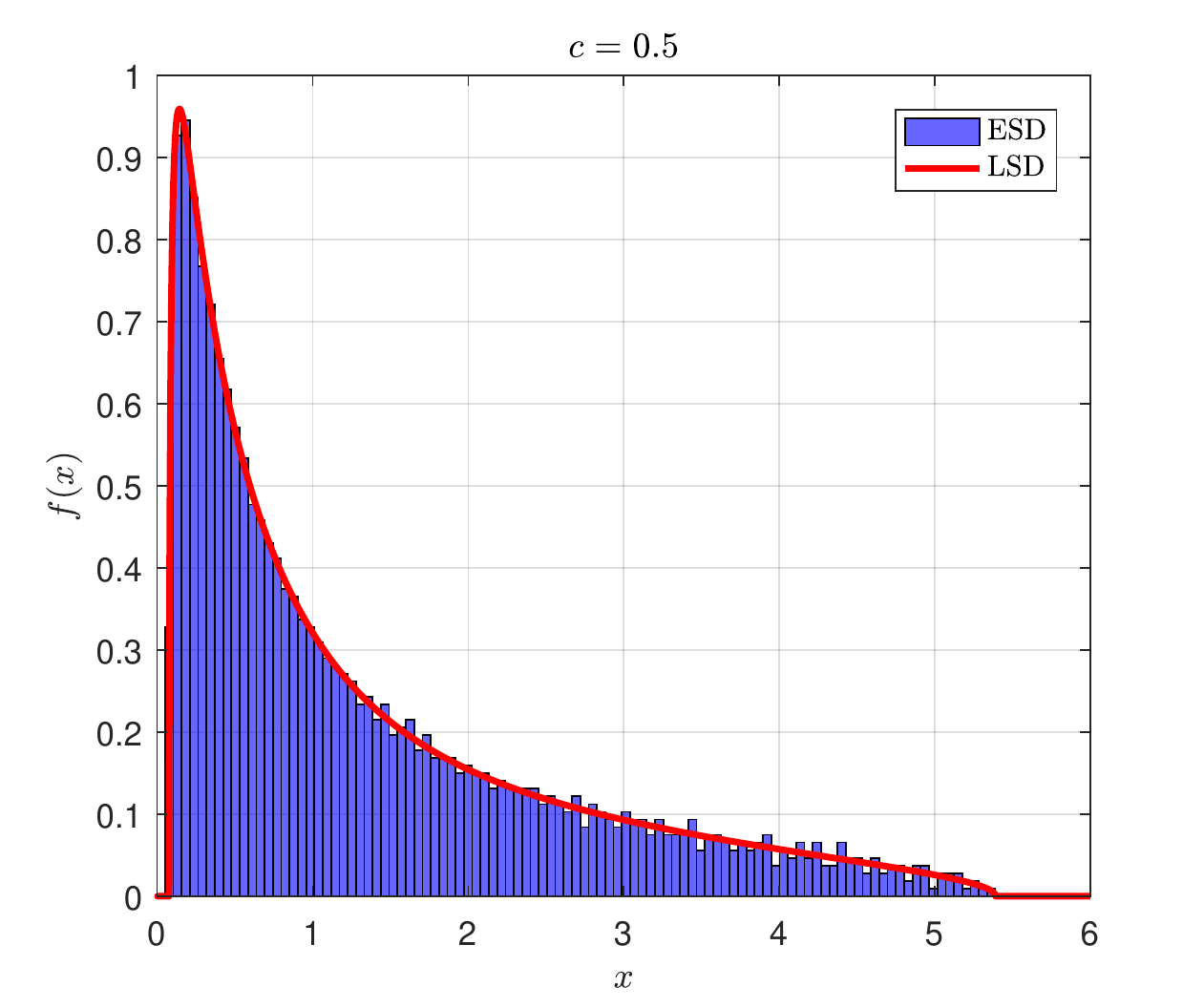}}
\caption{LSD and ESD of the heterogenous AR(1) process $U[i,t] = r_i U[i,t-1] + \xi[i,t]$ and the homogeneous AR(1) process $\bar{U}[i,t] = \bar{r} \bar{U} [i,t-1] + \eta[i,t]$.}
\label{F:LSD_ESD_AR_Residual}
\end{figure}

We demonstrate in Fig.~\ref{F:LSD_ESD_AR_Residual} the LSD and ESDs associated with the heterogenous and homogeneous AR(1) processes. The setting is the same as that in Fig.~\ref{F:MeanFieldModel}. We find that there both exists a good agreement between the theoretical LSD and ESD.

\subsection{Low-rank spiked sample covariance matrices}
We now look more closely at the eigenvalue phase transition associated with a low-rank spiked sample covariance model $Y_N = X_N + C_N$, where $C_N$ is the sample covariance matrix of a separable spatio-temporal process, i.e.,
\[
C_N = \frac{1}{T} \left( \Sigma^{s}_{N} \right)^{\frac{1}{2}} W_{N,T} \Sigma^{t}_{T} W_{N,T}^{\text{T}} \left( \Sigma^{s}_{N} \right)^{\frac{1}{2}} .
\]
We make the same assumptions as before: $X_N$ is with rank $r \ll N$ and non-negative eigenvalues $\theta_1 \ge \ldots \ge \theta_r \ge 0$, $W_{N \times T}$ is a Gaussian random matrix, $\Sigma^{s}_{N} = \alpha I_N + \beta V_{N}$, $\Sigma^{t}_{T} [a, b] = r^{ \left| {a-b} \right |}$ and $C_N \xrightarrow{distr} C$.

We first study the phase transition phenomenon for observation. Theorem~\ref{T:EigenvalueNonlinearMappingForward} states that, for any $\theta_i > \frac{1}{ \mathcal{G}_{\mu_C}(a^+)}$, there exists a corresponding eigenvalue of $Y_N$ such that $\lambda_i (Y_N) \xrightarrow{\text{a.s.}} \mathcal{G}_{\mu_C}^{-1} (\frac{1}{\theta_i})$. In fact, with the inverse relationship  $\mathcal{G}_{C} ^{-1} (\mathcal{G}_{C} (z)) = z$, we have the equation that $\mathcal{G}_{C}^{-1} (z)$ satisfies
\[
\begin{aligned}
1 = & \frac{c z}{2} \left( \gamma + \sqrt{ \gamma^2 - 1 + \frac{1}{c^2 \left( \mathcal{G}_{C}^{-1} (z) z - 1 \right)^2} }  \right) \\
& \left( \alpha + \sqrt{\alpha^2 + 4 \beta^2 \left( \mathcal{G}_{C}^{-1} (z) z - 1 \right) } \right) .
\end{aligned}
\]
With some algebras, it could be written as a 4-th order polynomial
\begin{equation}\label{E:G_Inverse_4order}
\sum_{i = 0}^{4} q_i \mathcal{G}_{C}^{-1} (z) ^i = 0 ,
\end{equation}
where
\begin{equation*}
\begin{aligned}
q_4 = & \beta^4 z^8 , \\
q_3 = & \frac{ 2 \beta^2 z^5 }{ c^2 } - \frac{4 \beta^2 \gamma^2 z^5}{c^2} + \frac{2 \alpha \beta^2 \gamma z^6}{c} - 4 \beta^4 z^7 , \\
q_2 = & \frac{z^2}{c^4} - \frac{2 \alpha \gamma z^3}{c^3} + \frac{\alpha^2 z^4}{c^2} - \frac{ 6 \beta^2 z^4 }{c^2} + \frac{12 \beta^2 \gamma^2 z^4}{c^2}  \\
& - \frac{6 \alpha \beta^2 \gamma z^5}{c} + 6 \beta^4 z^6 - \frac{2 \beta^4 z^6}{c^2} , \\
q_1 = & -\frac{2 z}{c^4} + \frac{4 \alpha \gamma z^2}{c^3} - \frac{2 \beta^2 z^3}{c^4} - \frac{2 \alpha^2 z^3}{c^2} + \frac{6 \beta^2 z^3}{c^2} \\
& - \frac{12 \beta^2 \gamma^2 z^3}{c^2} - \frac{2 \alpha \beta^2 \gamma z^4}{c^3} + \frac{6 \alpha \beta^2 \gamma z^4}{c} - 4 \beta^4 z^5 + \frac{4 \beta^4 z^5}{c^2} ,\\
q_0 = & \frac{1}{c^4} - \frac{2 \alpha \gamma z}{c^3} - \frac{\alpha^2 z^2}{c^4} + \frac{2 \beta^2 z^2}{c^4} + \frac{\alpha^2 z^2}{c^2} - \frac{2 \beta^2 z^2}{c^2} + \frac{4 \beta^2 \gamma^2 z^2}{c^2} \\
& + \frac{2 \alpha \beta^2 \gamma z^3}{c^3} - \frac{2 \alpha \beta^2 \gamma z^3}{c} + \beta^4 z^4 + \frac{\beta^4 z^4}{c^4} - \frac{2 \beta^4 z^4}{c^2} .
\end{aligned}
\end{equation*}
For general parameters, Equation~\ref{E:G_Inverse_4order} cannot be solved analytically. In much the same manner as solving Equation~\ref{E:G_Order8b}, we determine the roots of Equation~\ref{E:G_Inverse_4order} numerically.

\begin{figure}
\centering
\subfigure[Upward bias in the large eigenvalues.]{
\label{F:NonlinearMappingAR1:1}
\includegraphics[width=0.45\columnwidth]{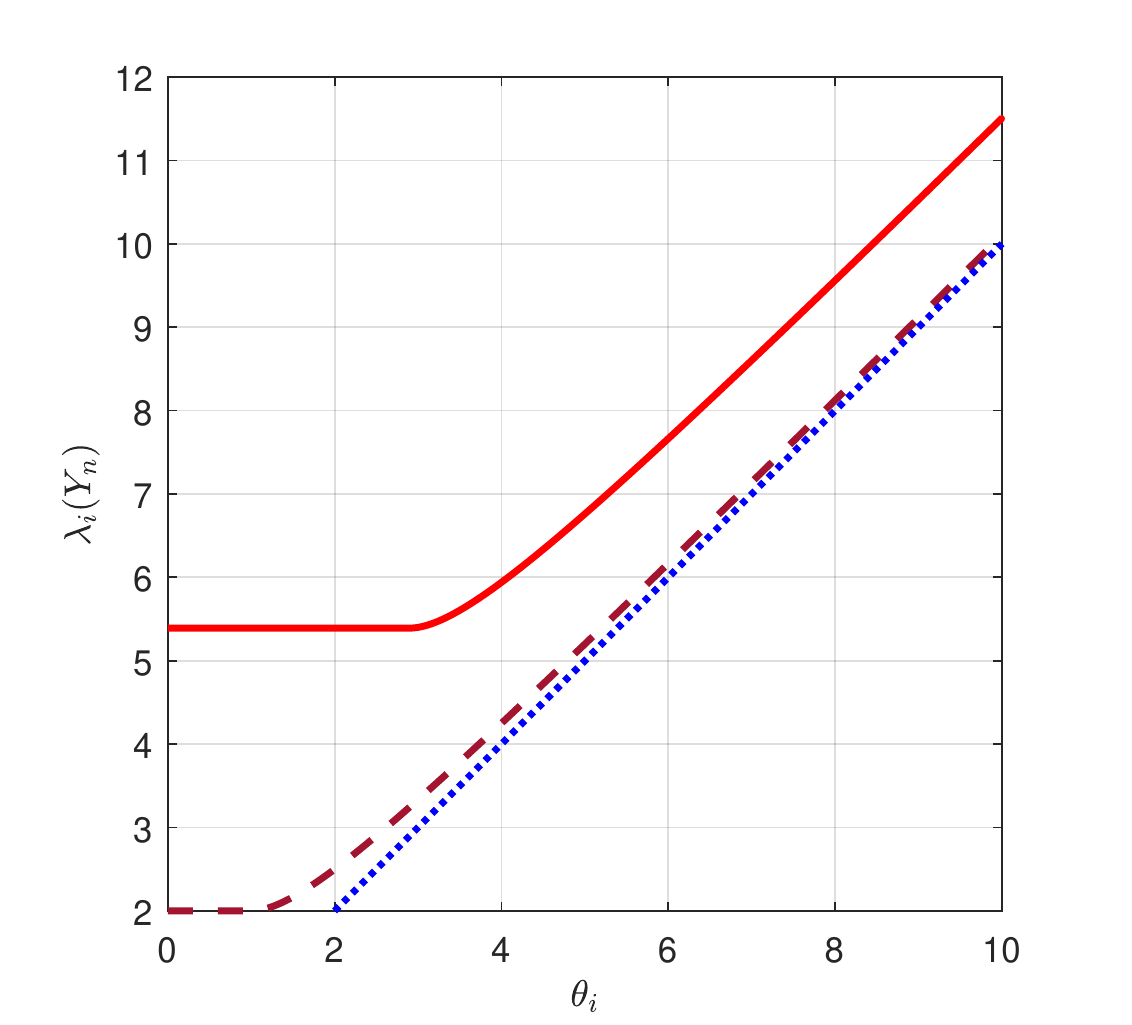}}
\subfigure[Nonlinear shrinkage estimator of the large eigenvalues.]{
\label{F:NonlinearMappingAR1:2}
\includegraphics[width=0.45\columnwidth]{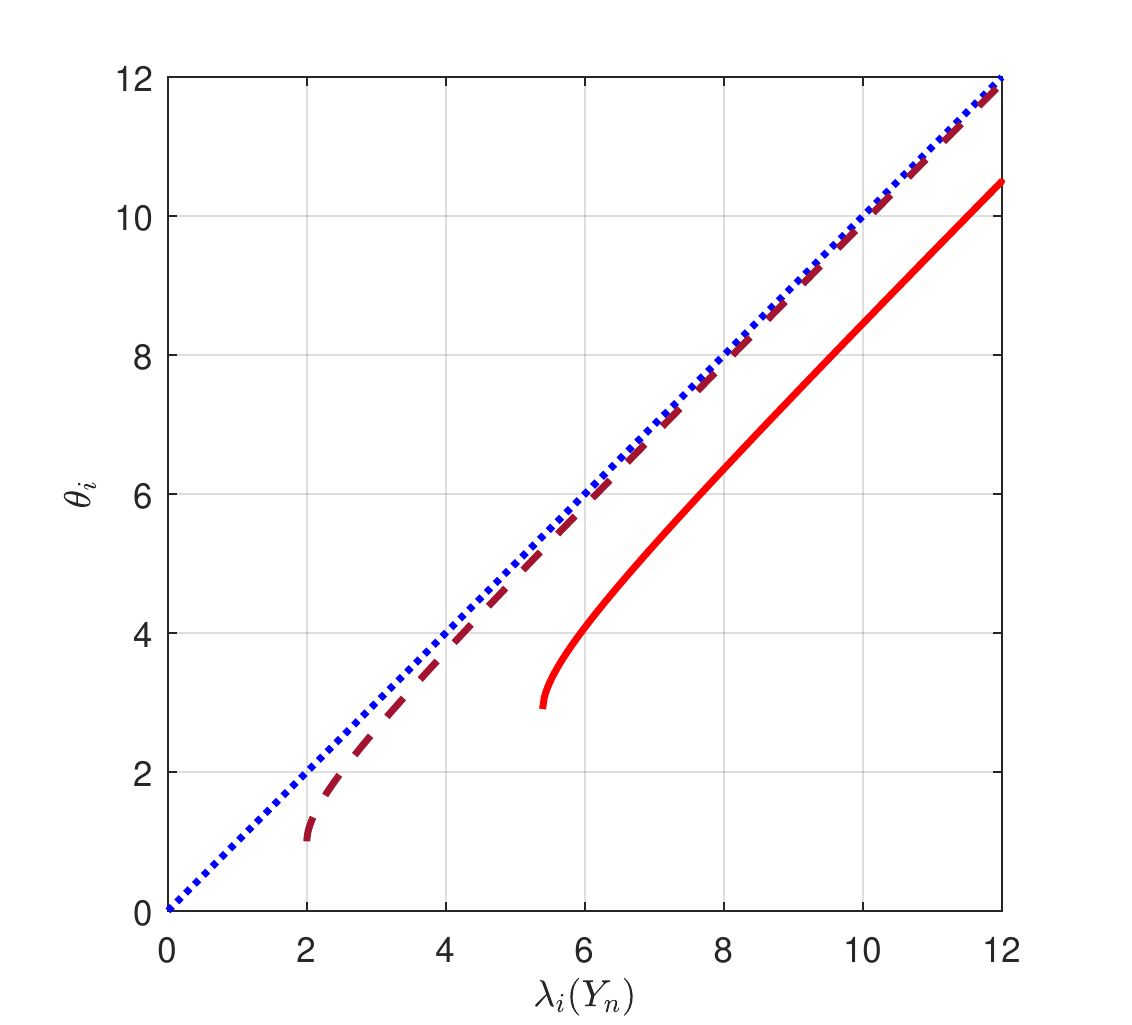}}
\caption{Eigenvalue phase transitions for low-rank spiked sample covariance model. The dashed brown is the mappings for noise covariance matrix being a Wigner matrix for easy comparison.}
\label{F:NonlinearMappingAR1}
\end{figure}

The phase transition phenomena are demonstrated in Fig.~\ref{F:NonlinearMappingAR1}, from which we can see there exists a greater amount of upward bias for observation (Fig.~\ref{F:NonlinearMappingAR1:1}) and shrinkage for estimation (Fig.~\ref{F:NonlinearMappingAR1:2}). It is also worth noting that the uncertain region for eigenvalue estimation becomes larger (Fig.~\ref{F:NonlinearMappingAR1:2}). This is primarily because the eigenvalues of sample covariance matrix of the form $C_N = \frac{1}{T} \left( \Sigma^{s}_{N} \right)^{\frac{1}{2}} W_{N,T} \Sigma^{t}_{T} W_{N,T}^{\text{T}} \left( \Sigma^{s}_{N} \right)^{\frac{1}{2}}$ are more spread out than those of Wigner matrix (comparison of Fig.~\ref{F:LSD_ESD_AR_Residual} and Theorem~\ref{T:SemiCircleLaw}).

At last, we will have a very brief discussion on the performance of nonlinear shrinkage estimators for top eigenvalues. As indicated in Section~\ref{SS:LowRankWigner}, for $N$ sufficiently large, if $\lambda_i (Y_N) >  \sup \left( \text{supp} ( \mu_{C} ) \right)$, $\widehat{\theta_{i}} = \frac{1}{\mathcal{G}_{\mu_C} (\lambda_i (Y_N))}$ is a nonlinear estimator for the $i$-th largest eigenvalue of $X_N$. Since we have no explicit expression for $\mathcal{G}_{\mu_C}(\cdot)$, the performance of the estimator is verified by  numerical experiments.

\begin{figure}
\centering
\subfigure[Heterogenous process]{
\label{F:EstimatedEigenvalue_AR1:1}
\includegraphics[width=0.47\columnwidth]{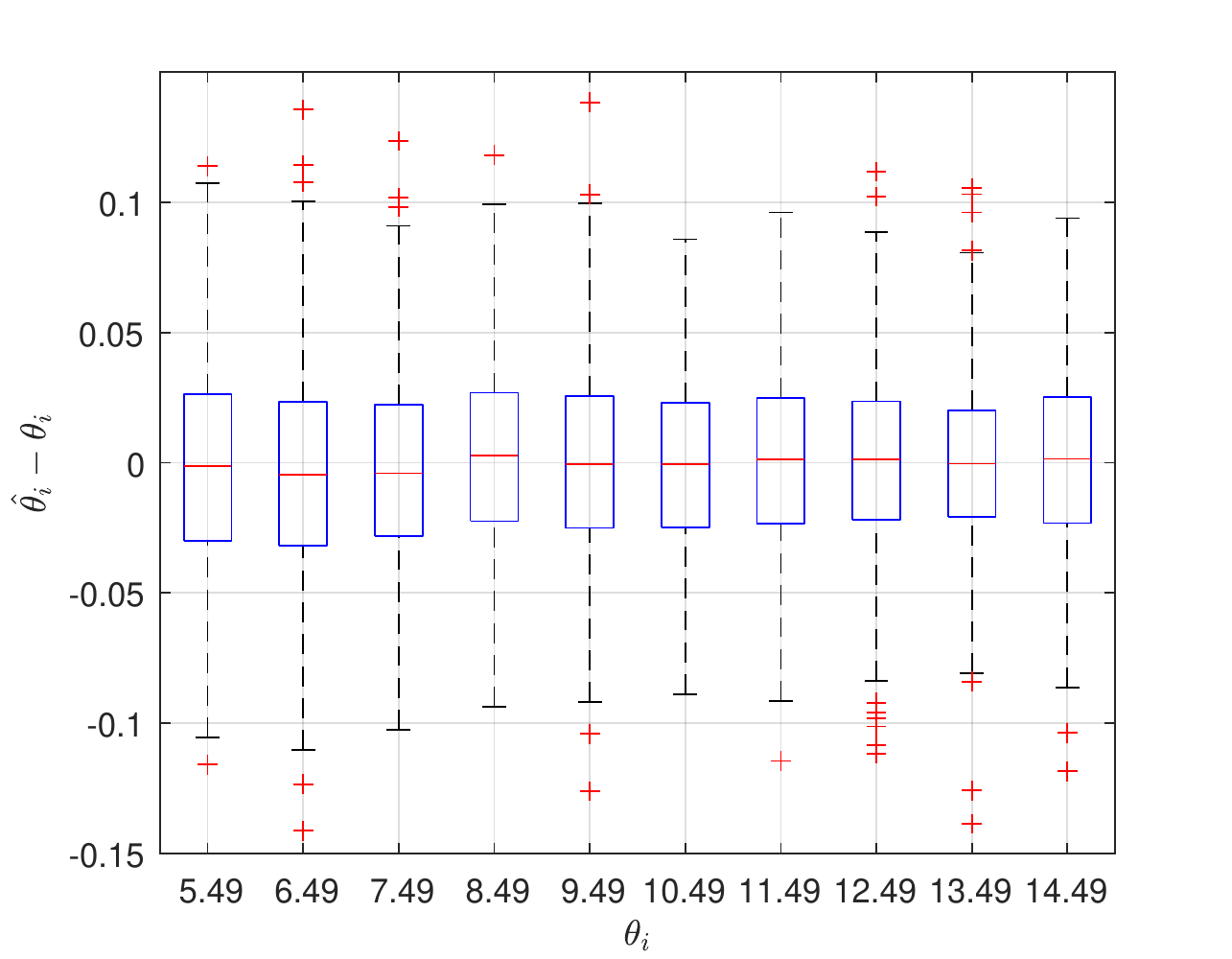}}
\subfigure[Homogeneous process]{
\label{F:EstimatedEigenvalue_AR1:2}
\includegraphics[width=0.47\columnwidth]{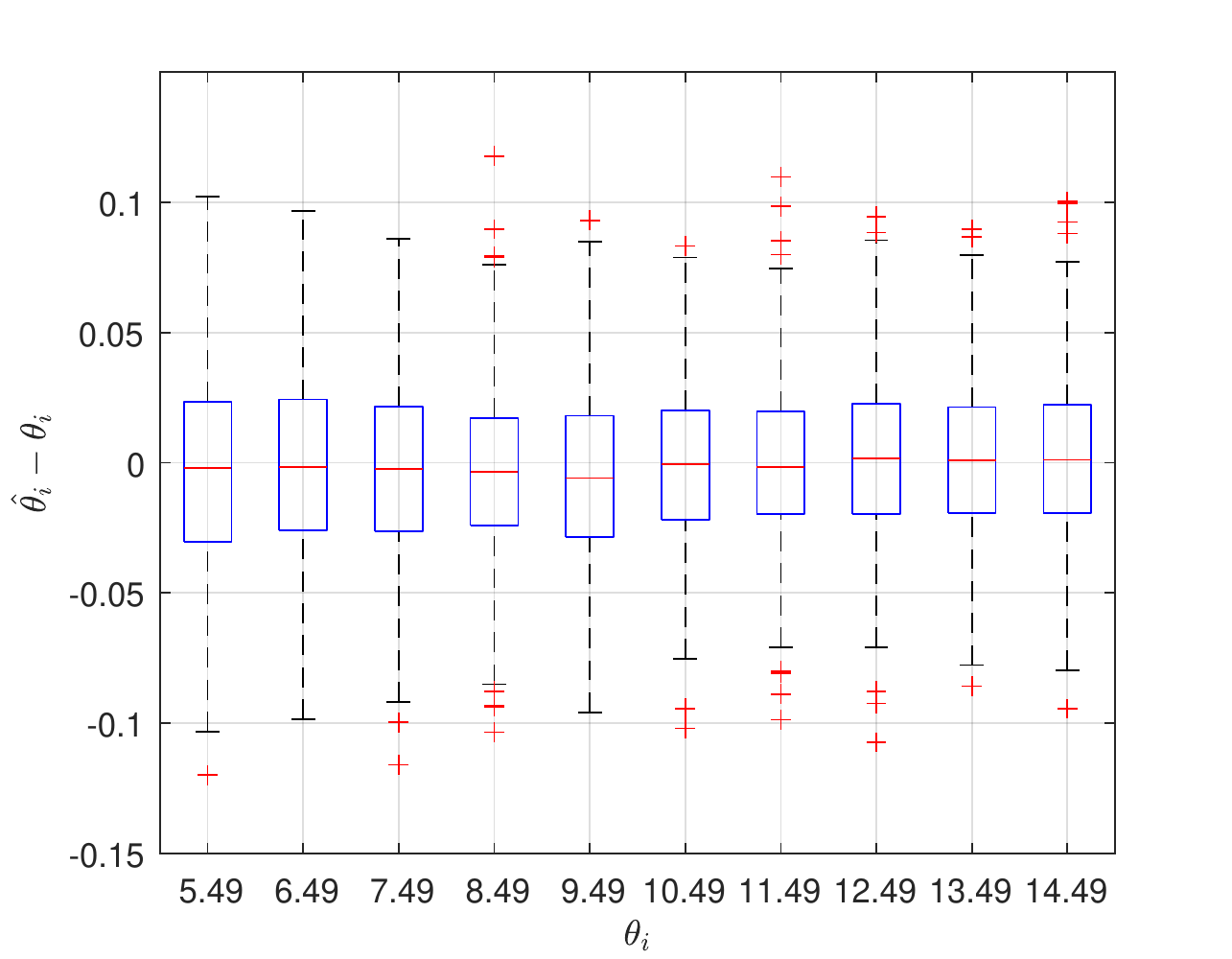}}
\caption{Boxplots of the error of the estimator $\widehat{\theta_{i}} = \frac{1}{\mathcal{G}_{\mu_C} (\lambda_i (Y_N))}$.}
\label{F:EstimatedEigenvalue_AR1}
\end{figure}

Fig.~\ref{F:EstimatedEigenvalue_AR1} shows boxplots of the error of the estimator $\widehat{\theta_{i}} = \frac{1}{\mathcal{G}_{\mu_C} (\lambda_i (Y_N))}$ for $C_N$ being sample covariance matrices associated with heterogenous (Fig.~\ref{F:EstimatedEigenvalue_AR1:1}) and homogeneous (Fig.~\ref{F:EstimatedEigenvalue_AR1:2}) AR(1) processes (N = 4000, T = 8000), respectivley. We let $X_N$ be a deterministic matrix with rank 10. For each $\theta_i$, a boxplot based on 500 trials is generated. This procedure results in 10 boxplots that are displayed side-by-side. We can conclude from these boxplots the nonlinear shrinkage estimator is pretty accurate.

%

\appendices

\ifCLASSOPTIONcaptionsoff
  \newpage
\fi



%
\bibliographystyle{IEEEtran}
\bibliography{References}

%

%





\end{document}